\documentstyle{amsart}
\newtheorem{thm}{Theorem}[section]
\newtheorem{lemma}{Lemma}[section]
\newtheorem{prop}{Proposition}[section]
\newtheorem{cor}{Corollary}[section]

\theoremstyle{definition}

\theoremstyle{remark}

\newtheorem{remark}{Remark}

\def \mrn {{\mathbb R}^n}
\def \msn {{\mathbb S}^n}
\def \mr {{\mathbb R}}
\def \ms {{\mathbb S}}
\def \mcf {{\mathcal F}}
\def \mc {{\mathbb C}}
\def \mh {{\mathbb H}}
\def \mn {{\mathbb N}}
\def \mz {{\mathbb Z}}

\def \eps {\epsilon}   
\def \la {\lambda}   
\def \lan {\langle}   
\def \ran {\rangle}   
\def \del {\delta}   

\newcommand{\Id}{\operatorname{Id}}

\def \ha {\frac{1}{2}}
\def \oq {\frac{1}{4}}

\def \p {\partial}

\def \sub {\subset}

\def \rao#1 {\frac{\p}{\p #1} #1}

\def \fdl#1#2{\fracwithdelims(){#1}{#2}}

\def \laphd {\tilde{\Delta}}

\numberwithin{equation}{section}
\title{Inverse Scattering on Asymptotically Hyperbolic Manifolds}
\author[Mark S. Joshi]{Mark S. Joshi}
\address{  \newline Department of Pure Mathematics and Mathematical
Statistics, \newline
University of Cambridge \newline
16 Mill Lane, \newline 
Cambridge CB2 1SB, England, U.K.}
\email{joshi@@dpmms.cam.ac.uk}
\author[Ant\^{o}nio S\'{a} Barreto]{Ant\^{o}nio S\'{a} Barreto}
\address{\newline Department of Mathematics  \newline
Purdue University, \newline
West Lafayette IN 47907, Indiana, U.S.A.}
\email{sabarre@@math.purdue.edu}
\begin{document}
\begin{abstract}
Scattering is defined on compact manifolds with boundary which are 
equipped  with an asymptotically
hyperbolic metric, $g.$ A model form is established for such metrics
close to the boundary.
It is shown that the scattering matrix at
energy $\zeta$ exists and is a pseudo-differential operator of order
$2\zeta+1 - \dim X.$ 
The symbol of the scattering matrix is then used
to show that except for a countable set of energies 
the scattering matrix at one energy determines the diffeomorphism class
of the metric modulo terms vanishing to infinite order at the
boundary. An analogous result is proved for potential scattering.
The total symbol is computed  when the manifold is
hyperbolic or is of product type modulo terms vanishing to infinite order 
at the boundary.  The same methods are then applied to studying
inverse scattering on the Schwarzschild and De Sitter-Schwarzschild models of
black holes.
\end{abstract}
\maketitle
\section{Introduction}

In this paper, we study scattering for Schr\"odinger operators
on asymptotically hyperbolic
manifolds. In particular, we show that the scattering matrix at energy
$\zeta$ is a
pseudo-differential operator of order $2 \Re \zeta -n$  (really complex 
order $2\zeta-n,$) where the dimension of the manifold is $n+1.$ 
We then show that
the total symbol of this operator is
determined locally by the metric and the potential and that,
except for
a discrete set of energies, the asymptotics of either the metric 
or the 
potential can be recovered from the scattering matrix at one energy.
 This also allows us to characterize the total symbol of
the scattering matrix in the case where the manifold is actually
hyperbolic, or when it is almost of product type.

We remark that the fact that the scattering matrix at energy
$\zeta$ is a pseudo-differential operator is a 
known result, see for example 
section 8.4 of \cite{mgst}.
However the proof in the general case does not seem to be available in the 
literature. The proof of several particular cases have been given,
see for example \cite{His}, \cite{guillop}, \cite{gz1} and \cite{per}
and references given there. The case $\Re \zeta=\frac{n}{2}$
is done in \cite{Bo}.

Recall that a compact manifold with boundary $(X, \p X)$ is asymptotically
hyperbolic if it can be equipped with a metric of the form
\begin{equation}
g = \frac{dx^2 + h(x,y,dx,dy)}{x^2}, \label{Intro1}
\end{equation} 
where $h_{|x=0},$  is independent of $dx$ for some boundary defining
function $x,$ and a 
product decomposition $X \sim \p X \times [0,\eps)$ near the boundary.
As observed in \cite{mm} this implies that along a smooth curve in
$X\setminus \p X,$ approaching a point $p \in \p X,$  the sectional
curvatures of $g$ approach $-1.$
We note that this form is invariant under multiplying $x$ by a
function of $y$ so there is no canonical metric on $\p X$ induced by
$g,$ but there is a natural conformal structure. The simplest example
of such a manifold is the hyperbolic space, ${\mh}^{n+1}$ and its quotients
by certain discrete group actions. 

Let $\Delta$ be the Laplacian on $X$ induced by $g.$  It will be shown in
section \ref{scatmat}
that  given a
function $f \in C^{\infty}(\p X)$ and $2\zeta \in \mc \setminus {\Bbb Z} ,$ 
$\zeta \in \mc \setminus (-\infty,0],$
which is not a pole 
of the meromorphic continuation of the resolvent,
there exists a unique solution of the
equation, $(\Delta +\zeta(\zeta-n))u=0$ of the form, 
\begin{equation}
u=x^{\zeta} f_{+} + x^{n-\zeta} f_{-},\label{sca1}
\end{equation} 
with $f_{+} , f_{-} \in C^{\infty}( X),$  and $f = f_{-}|_{\p x}.$ 
This is implicit in  \cite{maz1}, \cite{maz2} 
and is stated without a proof in \cite{mgst}.  A related result is also stated
in the introduction of \cite{mm}. The first terms of the expansion with
$\Re\zeta=\frac{n}{2}$ have been established in \cite{Bo}.

It is then natural to  define, for these values of $\zeta,$
 the scattering matrix to be  the map,
\begin{gather}
T(\zeta) : f \longmapsto {f_{+}}_{|\p X}.\label{scama}
\end{gather}
However the scattering matrix is then (mildly) dependent on the choice of 
boundary defining function $x$ and so we instead define it as a map,
\begin{equation}
S(\zeta) : C^{\infty}(\p X, |N^{*}(\p X)|^{n-\zeta}) \longrightarrow  
C^{\infty}(\p X, |N^{*}(\p X)|^{\zeta}) \label{sscatm}
\end{equation}
via $S(\zeta)( f |dx|^{n-\zeta} ) = (T(\zeta) f) |dx|^{\zeta}$ and it is then 
invariant. 
The same statements hold if we add  a short range potential to the
Laplacian  that is,
in this context, 
a real-valued function which is  smooth up to the boundary and
vanishes there. Whilst this definition can not make sense for $\zeta$
such that $2\zeta \in {\Bbb Z},$ as the decomposition \eqref{sca1} can
not then be unique, and the uniqueness of the expansion \eqref{sca1} is
not established in section \ref{scatmat} for $\zeta \in (-\infty,0],$ 
we shall see in Proposition \ref{PO7}, that the scattering matrix can be
defined as a restriction of the resolvent. This allows a meromorphic
continuation of $S(\zeta)$ across points which are not poles of \eqref{PO8}.

For simplicity, we work in a product decomposition such that 
\begin{gather}
 g = \frac{dx^2 + h(x,y,dy) }{x^2}.\label{pdec}
\end{gather}
The existence of such a model form is established in section 2. This
yields a trivialization of the normal bundle which we work with. 
In section \ref{thps} we prove 
\begin{thm}\label{pssm0} Let $(X, \p X)$ be a smooth manifold with
boundary. Suppose $g$ induces an asymptotically hyperbolic structure
on $X$ and that
$g = \frac{dx^2 + h(x,y,dy) }{x^2},$  with respect to some product 
decomposition near $\p X$.  Let $V \in C^{\infty}(X)$ be a short range 
potential and let $\zeta\in \mc$ be such that 
the scattering matrix,$S(\zeta),$  associated to
$\Delta_{g} + V + \zeta(\zeta-n)$ is defined. Then 
$S(\zeta)  \in \Psi DO^{2\Re \zeta -n},$ and its principal symbol equals
$C(\zeta) |\xi|^{2\zeta-n},$
where $|\xi|$ is the length of the co-vector $\xi$ induced by $h_0=h(0,y,dy)$
and $C(\zeta) = 2^{n-2\zeta} \frac{\Gamma(\frac{n}{2}
-\zeta)}{\Gamma(\zeta-\frac{n}{2})}.$
\end{thm}
This result has been established in \cite{Bo} for $\Re\zeta=\frac{n}{2}.$
As a direct consequence we obtain
\begin{cor}\label{pssm02} Let $(X, \p X)$ be a smooth manifold with boundary 
and let $p \in \p X.$ Suppose $g_1, g_2$
induce asymptotically hyperbolic structure on $X$ and that
$g_i = \frac{dx^2 + h_i(x,y,dy) }{x^2},$ $i=1,2,$ with respect 
to some product decomposition.  Let $S_{i}(\zeta)$ be the 
scattering matrix associated to $\Delta_{g_i} + V_i + \zeta(\zeta-n).$ 
There exists a discrete set $Q\subset \mc$ such that 
$S_1(\zeta)-S_2(\zeta) \in \Psi DO^{2\Re \zeta -n-1}$ for 
$\zeta \in \mc\setminus Q$ if and only if $h_1(0,y,dy)=h_2(0,y,dy).$
\end{cor}

We then analyze the difference of the scattering matrices when the metrics
$g_1$ and $g_2$ agree to order $k$ at the boundary.  We also prove in section 
\ref{thps}

\begin{thm}\label{pssm}
Let $(X, \p X)$ be a smooth manifold with boundary and let $p \in \p
X.$ Suppose $g_1, g_2$
induce asymptotically hyperbolic structures on $X$ and that
$g_i = \frac{dx^2 + h_i(x,y,dy) }{x^2},$ $i=1,2,$ with respect 
to some product decomposition. Moreover suppose that
$h_1 -h_2 = x^{k} L(y,dy)+ O(x^{k+1}),$ $k \geq 1,$ near $p$ and that 
$V_1, V_2$ are smooth
short range potentials such that
$V_1 - V_2 = x^{k}W(y) + O(x^{k+1}),$ near $p.$  Let $S_{i}(\zeta)$ be the 
scattering matrix
associated to $\Delta_{g_i} + V_i + \zeta(\zeta-n).$ 
We then have that, near $p,$
\begin{equation}
S_{1}(\zeta) - S_{2}(\zeta) \in \Psi DO^{2\Re \zeta -n -k},
\end{equation}
and the principal symbol of $S_{1}(\zeta) - S_{2}(\zeta)$ equals
\begin{equation}
A_1(k,\zeta)\sum \limits_{i,j} H_{ij}
\xi_i \xi_j |\xi|^{2 \zeta-n-k-2} +
A_{2}(k,\zeta)\left( W -\frac{1}{4}k(n-k)T\right)|\xi|^{2 \zeta-n-k},
 \label{psf}
\end{equation}
where $H= h_{0}^{-1} L h_{0}^{-1}$ as matrices, 
$h_0 = {h_1}|_{x=0}={h_2}|_{x=0},$ $T=\text{trace}\left(h_{0}^{-1}L\right),$
$|\xi|$ is the length of the co-vector $\xi$ induced by $h_0,$
and $A_1, A_2$ are meromorphic functions of $\zeta,$ which
 for $2\Re\zeta\geq \text{max}(n-k+1,k+2),$ are given by
\begin{gather}
\begin{gathered}
A_1(k,\zeta)= -\pi^{\frac{n}{2}} 2^{k+2-2\zeta+n} 
\frac{\Gamma\left(\frac{k+2-2\zeta+n}{2}\right)}{\Gamma\left(-\frac{k+2-2\zeta}{2}\right)} \frac{C(\zeta)}{M(\zeta)}
T_1(k,\zeta), \\
A_2(k,\zeta)= \pi^{\frac{n}{2}} 2^{k-2\zeta+n}
\frac{\Gamma\left(\frac{k-2\zeta+n}{2}\right)}{\Gamma\left(-\frac{k-2\zeta}{2}\right)} \frac{C(\zeta)}{M(\zeta)}
T_2(k,\zeta),
\end{gathered} \label{psfnew}
\end{gather} 
where $C(\zeta)$ and $T_j(k,\zeta),$ $j=1,2$ are given by \eqref{sm1}, and
$M(\zeta)$ is given by Proposition \ref{G01}.
\end{thm}
As our construction shows that the singularities of the kernel of the
scattering matrix are determined locally, we deduce that for
hyperbolic manifolds the total symbol will agree with that for the
model hyperbolic space and we have:
\begin{thm} \label{product}
If $(X, \p X)$ is an infinite volume smooth hyperbolic manifold with
funnels, or if in some product decomposition the metric is a product modulo
terms vanishing to infinite order at the boundary
then  the scattering matrix is equal to 
$$ 2^{n-2\zeta} \frac{\Gamma(\frac{n}{2}
-\zeta)}{\Gamma(\zeta-\frac{n}{2})} \Delta_{\p X}^{\zeta-\frac{n}{2}},$$
 modulo smoothing. Here we have  chosen a defining function $x$ in
order to trivialize the normal bundle and to induce a metric on the
boundary, with respect to which we take $\Delta_{\p X}.$ 
\end{thm}

In the hyperbolic case, this result is due to Perry
\cite{per} (Perry's definition of the scattering matrix
was slightly different which caused an extra factor to be present.)
 We prove the result for almost product type
structures in
section \ref{sec:prod}. 

As consequences of Theorem \ref{pssm} we  have the inverse results:
\begin{cor}\label{cor1} Let $(X, \p X), g_{j}, S_{j}$ be as in 
Theorem \ref{pssm}, let 
$p\in \p X$ and
suppose that $V_1= V_2$ near $p.$ There
exists a discrete set $Q \sub \mc,$ such that if $\zeta \in \mc \setminus Q,$
and $S_{1}(\zeta) - S_{2}(\zeta) \in \Psi DO^{2\Re\zeta-n-k},$ $k\geq 1,$
 near $p,$ then
there exists a diffeomorphism $\psi$ of a neighbourhood $U \subset X$ of
$p,$ fixing 
$\p X,$ such that 
$\psi^{*} g_1 - g_2 = O(x^{k}).$ 
\end{cor}

\begin{cor}\label{cor2}
Let $(X, \p X), V_{j}, S_{j}$ be as in Theorem \ref{pssm}, let $p \in \p X$ and suppose that
$g_1= g_2$ near $p.$
There exists a discrete set $Q \sub \mc$ such that if 
$\zeta \in \mc\setminus Q$
and $S_{1}(\zeta) - S_{2}(\zeta) \in \Psi DO^{2\Re\zeta-n-k},$ $k\geq 0,$
 near $p,$ then $V_1 - V_2 =O(x^{k})$ near $p.$
\end{cor}
Of course intersecting over all $k,$ we see that off a countable set
of energies a metric or potential can be recovered modulo terms
vanishing to infinite order at the boundary. We will prove these Corollaries in section \ref{thps}, after the proof of
Theorems \ref{pssm0} and \ref{pssm}.

In  section \ref{schw} we give an application of these results, or rather of 
the methods used to prove them, to inverse scattering on the Schwarzschild 
and De Sitter-Schwarzschild model
of black holes.  We show that the Taylor series at the boundary of certain 
time independent perturbations of these models can be recovered from the
scattering matrix at a fixed energy.

Our approach is heavily influenced by the work of Guillop\'e and
Zworski, \cite{guillop}, \cite{gz1} and \cite{gz2}. In particular, we
compute the scattering matrix as a boundary value of the resolvent. To
do this we use the calculus developed by Mazzeo and Melrose
\cite{mm} of zero pseudo-differential operators in order 
to construct the resolvent. 

As in our work on asymptotically Euclidean scattering,
\cite{coulomb,smaginv,metric},  
 a key part of our approach is
to consider the principal symbol of the difference of the scattering
matrices rather than the lower order terms of the symbol of a single
operator, which allows
us to proceed more invariantly. We remark that whilst our results are
quite similar to those in the Euclidean case, the proofs and
underlying ideas are very different. The fundamental reason being that
in the asymptotically Euclidean category, as observed by 
Melrose \cite{masy}  and by Melrose and Zworski
\cite{smagai}, there is  propagation of
growth at infinity whilst  this does not occur in the
asymptotically hyperbolic category. This is reflected in the fact,
 proved in \cite{smagai}, that in the asymptotically Euclidean case
the scattering matrix is a
Fourier integral operator associated to the geodesic flow at time 
$\pi,$ whilst in the asymptotically hyperbolic manifold case it is a 
pseudo-differential operator and in the fact that the principal symbol of the
difference of the scattering matrices is locally determined by the
perturbation. See \cite{optics} for a discussion of a general
framework including both cases.  

There is a long history of scattering theory on hyperbolic manifolds
arising from the observation that the Eisenstein series for a Fuchsian
group is a generalized eigenfunction for the Laplacian on the
associated quotient of hyperbolic space - the fundamental reference
for this is \cite{lp1} where the finite volume case is studied. 
The study of the infinite volume case was initiated by Patterson in
\cite{patt}. There has been a wealth of results in both cases and we
refer the reader to \cite{gz1} for a comprehensive bibliography and to 
\cite{His} and \cite{mgst} for a review of the subject. There has been less 
work on
asymptotically hyperbolic spaces.  Mazzeo-Melrose, \cite{mm}, and
Mazzeo  \cite{maz1}, \cite{maz2} studied properties of the Laplacian
on such manifolds from which the properties of the scattering matrix 
proved in section 4 are implicit. In \cite{Bo} Borthwick showed
the continuous dependence of the scattering matrix on the metric.
 Agmon has also studied related questions, see \cite{ag1}, \cite{ag2}. 
Andersson, Chrusciel and Friedrichs
have studied solutions of the Einstein equations and related problems
on asymptotically hyperbolic spaces, \cite{and1}, \cite{and2}. There
appears to be no results in the literature on the inverse scattering
problem on asymptotically hyperbolic manifolds. Perry, \cite{per2},
has shown that for hyperbolic quotients in three dimensions by convex,
co-compact, torsion-free Kleinian groups with non-empty regular set,
that the scattering matrix determines the manifold. Borthwick, McRae
and Taylor have proved an associated rigidity result, \cite{borth1}.

We would like to thank Maciej Zworski for  explanations of 
hyperbolic scattering and helpful comments. We would also like to
thank Richard Melrose and Rafe Mazzeo for helpful conversations. We
are also grateful to Tanya Christiansen for explaining her computation
in the almost product case in the asymptotically Euclidean setting. 
This
work was initiated whilst visiting the Fields Institute and we would
like to thank that institution for its hospitality. This research was
partially supported by an EPSRC visiting fellowship. The second author
was also partly supported by NSF under grant DMS-9623175.

\section{A Model Form}

In this section, we establish a model form for asymptotically
hyperbolic metrics near infinity (the boundary.) This is very similar
in statement and proof to the model form for scattering metrics proved
in \cite{metric}. In the case where all sectional curvatures are equal to
$-1$ near the boundary, such normal form has been established in \cite{gz3}.

\begin{prop} \label{prop:modelform}
Let $(X, \p X)$ be a smooth manifold with boundary $\p X.$ And
suppose $g$ is a metric on $X$  such that
$$ g = \frac{dx^{2} + h(x,y,dx,dy)}{x^2},$$
in some product decomposition near $\p X,$ where $x$ is a defining function of
$\p X,$ with $h_{|x=0}$ independent of $dx.$ 
Then there exists a product decomposition,
$(\bar{x}, \bar{y}),$ 
near $\p X$  such that
\begin{gather}
 g = \frac{{d\bar{x}}^2 + h(\bar{x},\bar{y},d\bar{y})}{\bar{x}^2}.
\label{Mo1}
\end{gather}
\end{prop}
\begin{pf} First we prove this result modulo terms that vanish to
infinite order at $x=0.$
It is enough to show the existence of a sequence of diffeomorphisms,
$\psi_k ,$  of $\p X \times [0,\eps) $ such that 
$$\psi_{k}^{*} g = \frac{{d\bar{x}}^2 + h(x,y,d\bar{y})}{\bar{x}^2} +
O(\bar{x}^{k}), $$ and 
$$f_{k} = \psi_{k-1}^{-1} \psi_{k},$$ fixes the boundary to order $k+1.$ 
This is enough as a diffeomorphism $\psi$ can then be picked, using the
Borel lemma, of which the $l^{th}$ term in the Taylor series will
agree with that of 
$\psi_k$ for $l \leq k,$ for all $l.$ 

Suppose $\psi_{k-1}$ has been constructed. We show how to pick
$f_k$ so that 
$$f_{k}^{*} \psi_{k-1}^{*} g = \frac{{d\bar{x}}^2 + h(x,y,d\bar{y})}{\bar{x}^2} +
O(\bar{x}^{k}).$$ Putting $\psi_k = \psi_{k-1} f_{k},$ our result
follows. 

We work in local coordinates on the boundary. We shall see that the choice
of the next term in the Taylor series is actually unique so there is
no problem patching these local computations together.
So suppose we have,
\begin{equation*}
\psi_{k-1}^{*} g = \frac{dx^2 + h(x,y,dy)}{x^2} + x^{k} (
\alpha(y) dx^2)  + x^{k} (\sum \beta_{j}(y) dx dy_{j}) +
O(x^{k+1}).
\end{equation*}
Putting
$$ x= \bar{x} + \gamma(\bar{y}) \bar{x}^{l},$$ and 
$$ y= \bar{y} + \delta(\bar{y}) \bar{x}^{l},$$ we have
$$ dx = d \bar{x} + l \bar{x}^{l-1} \gamma(\bar{y}) d\bar{x} +
\bar{x}^{l} \frac{\p
\gamma}{\p \bar{y}} d\bar{y},$$
and 
$$dy = d\bar{y} + l \bar{x}^{l-1} \delta(\bar{y}) d\bar{x} + x^{l}
\frac{ \p \delta}{\p \bar{y}}(\bar{y}). d\bar{y}.$$ 

Now if $h(0,y,dy) = \sum h_{ij}(y) 
dy_{i} dy_{j}$ and $l=k+3,$ we see that
the metric becomes, modulo $ O(\bar{x}^{k+1})$ terms,  
$$\frac{d\bar{x}^2 + h(\bar{x},\bar{y},d\bar{y})}{\bar{x}^2} + \bar{x}^{k} (
\alpha(\bar{y}) d\bar{x}^2 ) + \bar{x}^{k} (\sum \beta_{j}(\bar{y})
d\bar{x} d\bar{y}_{j}) + 2(k+3) \gamma(\bar{y}) \bar{x}^{k}d\bar{x}^2
 +2(k+3)x^k\sum
h_{ij} \delta_{i} d\bar{y}_j dx.$$
Picking $\gamma(\bar{y})=-\frac{\alpha(\bar{y})}{k+3},$ and
as the form $h_{ij}(\bar{y})$ is non-degenerate, there is a unique choice of
$\delta$  such that  $2(k+3)x^k\sum h_{ij} \delta_{i} =-\beta_j.$ This
kills the terms of order $k$ in $dx^2$ and
$dx dy.$ Equation \eqref{Mo1}, modulo $O({\bar{x}}^\infty)$ terms, 
follows.

Having achieved the modulo form modulo $x^{\infty},$ which is all that is
necessary for the rest of this paper, we now show this form can be improved to
remove this error. If 
$$ g = \frac{dx^{2} + h(x,y,dy) + O(x^{\infty}) }{x^2},$$
then the geodesic flow is generated by the Hamiltonian function
$$ \sigma = x^{2} \tau^{2} + x^{2} \sum \limits_{i,j} h^{ij}(x,y,\xi)
+ {\cal O}(x^{\infty}).$$ 
Now if we work in rescaled zero coordinates, that is let $\bar{\tau} =
x \tau, \bar{\xi} = x \xi,$ and leave $(x,y)$ fixed, the canonical 1-form
$\alpha=\tau dx+\xi\cdot dy$ is rescaled to
$${}^{0}\alpha= \bar{\tau}\frac{dx}{x}+ \bar{\xi}\cdot \frac{dy}{x}.$$

The $0$-Hamilton vector field of $g,$ $H_g,$ is defined by

$$d{}^{0}\alpha(\bullet, H_g)=dg.$$

We then find that, modulo $O(x^\infty)$ terms,

$$H_g= 2\bar{\tau}\left(x\frac{\p}{\p x} + 
\bar{\xi}\cdot\frac{\p}{\p \xi}\right)-\left(2h^{-1}
+x\frac{\p}{\p x}h^{-1}\right)\frac{\p}{\p \bar{\tau}}+ xH_{h^{-1}},$$
where $h^{-1}(x,y,\bar{\xi})=\sum_{i,j} h^{ij}(x,y){\bar{\xi}}_i{\bar{\xi}}_j,$
and 
$H_{h^{-1}}=\sum_{i} \left(\frac{\p}{\p {\bar{\xi}}_i}h^{-1}\frac{\p}{\p y_i}
-\frac{\p}{\p y_i}h^{-1}\frac{\p}{\p {\bar{\xi}}_i}\right).$
 
This is of the form,
$$ H_g=2\bar{\tau} \left( x\frac{\p}{\p x} +
 \bar{\xi}.\frac{\p}{\p \bar{\xi}}
\right) 
+ {\cal O}(x^2+|\bar{\xi}|^2),$$
where $|\bar{\xi}|^2= h^{-1}(x,y,\bar{\xi}).$

Now if we restrict to the cosphere bundle, $$\bar{\tau}^{2} +
h^{-1}(x,y,\bar{\xi}) =1,$$ which is invariant under the flow,
we can re-express $\bar{\tau}$ in terms of
$(x,y,\bar{\xi}),$ and near $\bar{\tau}=-1$ the vector field becomes, 
$$H_g= -2 \left( x \frac{\p}{\p x} + \bar{\xi}.\frac{\p}{\p \bar{\xi}}
\right) + {\cal O}({\bar{\xi}}^{2}+x^{2}).$$ 

This forms a sink at $(x,\xi) =0$ and thus, by Theorem 7 of \cite{sell},
there exist local
coordinates $(x',\xi'),$ equal to $(x,\xi)$ to second order at $(x,\xi)
=0,$ which reduce the vector field to the form 
$$ -2 \left(x' \frac{\p}{\p x'} + \bar{\xi}.\frac{\p}{\p \xi'}
\right) + {\cal O}({\xi'}^{2}+{x'}^{2})\frac{\p}{\p y}.$$ 
We therefore see that any integral curve starting close enough to
$(x,\xi) = (0,0)$ will converge to $(0,0).$ 

So in particular if we take a hypersurface $S_{\eps} = \{ x=
\eps\},$ then the geodesics starting on the unit normals pointing to the
boundary will converge to $(x,\xi) =(0,0).$ We also have that the $x'$
derivatives of these geodesics will be non-zero so they can be
reparametrized in terms of $x'.$ In $x'>0,$ we can put $\theta=
\xi'/x',$ and use $(x',\theta, y)$ as coordinates. The form of the
vector field means that the angular coordinate $\theta$ will be
constant on geodesics and we have,
$$ \frac{dy}{dx'} = {\cal O}({x'}^2 + x' \theta^2) = {\cal O}(x').$$
Now as the finite time solution of an ODE, the point $y$ on the
boundary, which is the limit of the geodesic, will vary smoothly with
the start point on $S_{\eps}.$ We also see that, as the change in $y$ is
the integral of the derivative along the curve, that the Jacobian of
the map from $S_{\eps}$ to the boundary will be invertible for $\eps$
sufficiently small. This will also be true for the map to
hypersurfaces $S_{\eps'}, \eps' < \eps.$ 

So if we now take geodesic normal coordinates to the hypersurface
$S_{\eps},$ then these give us a map $\chi: \p X \times \mr_{+} \to X,$
which extends to the compactification of each to give a map, which is
smooth up to the boundary, and is a diffeomorphism in a neighbourhood
of the boundary. Now these coordinates are of the form,
$$ds^{2} + h(s,y,dy),$$ so if we put $X= e^{-s}$ we get coordinates in
a neighbourhood of infinity such that the metric is of the form,
$$ \frac{dX^{2}}{X^2} + h(X,y,dy).$$ Note that the change of
coordinates gives the correct compactification at the boundary to
ensure there has been no change of smooth structure there. Thus we
know the transformed metric must be of the form 
$$\frac{{dX}^2 + l(X,y,dy,dX)}{X^2}$$ with $l$ smooth up to $X=0,$ so
we conclude that 
that $h(X,y,dy) = \frac{k(X,y,dy)}{X^2}$ with smooth
up to $X=0$ and are we are
done. 
\end{pf}
We remark that the construction of the Taylor series in the first part
of the proof
gave a boundary defining function of the form $\bar{x}= x + {\cal
O}(x^2)$ and that the rest of the Taylor series was then
determined; one could however start with a different defining function
$\alpha(y) x.$  This contrasts with the case of a scattering metric
where the $x^{1}$ term is fixed by the metric but the $x^{2}$ term can be
chosen.

\section{Constructing the Resolvents}\label{S.res}

In this section, we review the construction of the resolvent on an
asymptotically hyperbolic space due to Mazzeo and Melrose and show how
to modify it to obtain information about the difference of two
resolvents associated to data which agree to some order at the
boundary. Our account is necessarily brief and we concentrate on
explaining where our construction differs from theirs and refer the
reader to their paper \cite{mm} for further details. We shall work
with half-densities throughout as they give better invariance
properties. 

We recall that a Riemannian metric $g$ on a manifold $Y$ induces a
canonical trivialization of the $1-$ density bundle by taking
$\omega=\sqrt{\delta} |dy|$ where $\delta$ is the determinant of
$g_{ij}$ in the local
coordinates $y$ on $Y.$ The square root of this is then a natural
trivialization of the half-density bundle. We then have a natural
Laplacian, $\laphd,$  acting on half-densities by 
$$\laphd(f \omega^{\ha} )= \Delta(f) \omega^{\ha}.$$ 

Mazzeo and Melrose showed that the resolvent could be meromorphically
continued to the entire complex plane and that it could be constructed
in a certain class of ``zero'' pseudo-differential operators. A
``zero'' vector field is a vector vanishing at the boundary and a
``zero'' differential operator is a composition of such vector fields.
The most important example being the Laplacian associated to an
asymptotically hyperbolic metric. 

``Zero'' pseudo-differential 
operators have kernels living on the blown-up space $X \times_{0} X.$
This is the space  obtained by blowing up $X \times X$ along the
diagonal, $\Delta_{\p X},$ of $\p X
\times \p X.$ We recall that blow-up is really just an invariant way
of introducing polar coordinates and that a function is smooth on the
space $X \times_{0} X$ if it is smooth in polar coordinates about
$\Delta_{\p X} .$ As a set, $X \times_{0} X$ is $X \times X$ with 
$\Delta_{ \p X}$ replaced by the interior pointing portion of its normal
bundle. Let
$$ \beta :X \times_{0} X \longrightarrow X \times X $$
denote the blow-down map.
If $(x,y)$ are coordinates in a product decomposition of
$X$ near $\p X,$  and we let $(x',y')$ be the corresponding
coordinates on a second copy of $X,$ then $R = (x^2 + {x'}^2 +
(y-y')^2)^{\ha}$ is a defining function for the new face which we call
the front face. The functions $\rho = x/R$ and $\rho' = x'/R$ are then
defining functions for the other two boundary faces which we call the
top and bottom spaces respectively. One advantage of working on this
blown-up space is that the lift of the diagonal of $X\times X$ only meets the 
front 
face of the
blown-up space and is disjoint from the other two boundary faces.

To define the space of ``zero'' pseudo-differential operators, Mazzeo
and Melrose defined a bundle $\Gamma_{0}(X),$  whose sections are
smooth multiples of the Riemannian density.  Note that
for the Riemannian
structure \eqref{Intro1}, the natural density  is singular at $\p X.$
In local coordinates $(x,y),$ where $x$ is a defining function of the boundary,
it is given by 
\begin{gather*}
h(x,y) \frac{dx}{x}\frac{dy}{x^n}, \;\ h \in C^\infty(X), \;\ h \not =0.
\end{gather*}
We denote $\Gamma_{0}^{\ha}(X)$ the analogous bundle of half-densities.
Similarly we define the bundle $\Gamma_{0}^{\ha}(X \times X).$
The bundle $\Gamma_{0}^{\ha}$ over $X \times_{0} X$ is then defined to be the
lift of $\Gamma_{0}^{\ha}(X \times X)$ under the blow-down map.

A ``small''
zero-pseudo-differential operator of order $m,$  is then an operator on $X$ of which
the Schwartz kernel when lifted to $X \times_{0} X$ vanishes to
infinite order at the top and bottom faces, and is the restriction of
a section of $\Gamma_{0}^{\ha}$ over the double across the front face,
which is conormal to
the lifted diagonal of order $m.$ In the interior, these are of course
just the usual class of pseudo-differential operators acting on
half-densities. The space of
these kernels will be denoted $K_{0}^{m}(X)$ and the corresponding
operators by $\Psi_0^m(X,\Gamma_0(X)).$

The ``large class'' $\Psi_{0}^{m,s,t}(X), s,t \in \mc$ is then defined
to be operators 
which have Schwartz kernels that are equal to an
element of $K_{0}^{m}(X)$ plus a smooth function of the form
$\rho^{s}{\rho'}^{t} f$ with $f \in C^{\infty}(X \times_{0} X,
\Gamma_{0}^{\ha})$ and smooth up to the boundary. This space then has
three natural filtrations but it will also be important to consider a
fourth which is the order of vanishing at the front face, so we
commonly work with operators with kernels in the class $R^{k}
\Psi_{0}^{m,s,t}(X).$  In \cite{mm}, Mazzeo and Melrose show that the
resolvent of $\Delta + \zeta(\zeta -n)$ has a meromorphic extension to
all of $\mc$ 
and that it lies in $\Psi_{0}^{-2,\zeta,\zeta}(X).$ 

The ordinary symbol map expressing the lead singularity at the
diagonal, extends to this class and is a homogeneous section of the
zero-cotangent bundle - that is the dual bundle to the space of vector
fields vanishing at the boundary. There is also a second natural
symbol map which is called the normal operator. This is obtained by
restricting the Schwartz kernel to the front face and therefore
expresses the lead term there, which is therefore a section of the
bundle $\Gamma_{0}^{\ha}(X\times_0 X)$ restricted to that face. 

Let $p \in \p X$ and let $X_p$ be the inward pointing vectors in
$T_{p}(X).$  This a manifold with boundary and has a metric 
$$g_p = (dx)^{-2} h_p,$$ where $g= (dx)^{-2} h,$ making it isometric to
the hyperbolic upper half-plane. ( We regard $h_p$ and $dx$ as linear
functions on the tangent space $X_p .$ )
Mazzeo and Melrose observed that the leaf of the front face above a
point $p$ is naturally
isomorphic to $X_p,$ using a natural group action on the front
face. This group action is obtained by lifting the action of the
subgroup of the general linear group of the boundary of $X_p$ to the
normal bundle of $X_p,$ as a leaf of the front face is just a quarter
of the normal bundle over $p.$ 

It is also observed in \cite{mm}, that the restriction of
$\Gamma_{0}^{\ha}(X\times_0 X)$ to the front face is
canonically trivial, and then can act as a convolution operator
using the natural group structure on the front face.  As mentioned above,
the fibre of
the front face above a point $p$ can be identified with  $X_p.$ 
 If we take local
coordinates $(x,y)$ with $x$ a boundary defining function and denote
the natural corresponding linear coordinates on $X_p$ by $(x,y)$
also. Let $(x',y')$ be the same coordinates on the right factor in $X
\times X$ and let $s = x/x', z = (y-y')/x.$ Then if the Schwartz
kernel of a map $B$ is $k(x',y',s,z) \gamma$ with $\gamma = \left|
\frac{ds dz dx dy}{ sx^{n+1}} \right|^{\ha},$ the normal operator is
given at $p=(0,\bar{y})$ by
\begin{equation}
\left[ N_{p}(B)(f \mu) \right] = \int k(0,\bar{y},s,z) f\left( \frac{x}{s},
y - \frac{x}{s} z\right) \frac{ds}{s} dz .\mu, 
\end{equation}
where $d\mu = \left| \frac{dx}{x} \frac{dy}{x^n} \right|^{1/2}.$ 

 In fact, Mazzeo and Melrose only used
the normal operator for terms in $\Psi_{0}^{-\infty,s,t}(X)$ but it works
equally well for terms in $\Psi_{0}^{m,s,t}(X),$ see Theorem 4.16 of 
\cite{mm}. The main difference
being that the normal operator instead of being a smooth half-density 
on the front face,  now has a conormal singularity at the
centre, i.e the intersection of the lift of the diagonal of $X \times X$
with the front face.
The normal operator will of course have growth at the
boundaries of the front face according to $s,t.$ In particular it will
be in the space ${\cal A}^{s,t}$ of half-densities growing of order
$s$ at the top edge and of order $t$ at the bottom.

The important fact is that the normal operator of a zero differential
operator is obtained by freezing the coefficients at a point on the
boundary and the normal operator of the Laplacian is just the
Laplacian of the induced metric on the space $X_p.$ 
As a short range potential vanishes at the boundary, if
 $P(\zeta) = \Delta +V+ \zeta(\zeta-n),$ with $V$ short
range, and $Q \in \Psi_{0}^{m,s,t}(X),$ we thus have
that
\begin{equation}
N_{p}(P(\zeta) Q) = (\Delta_{p}+\zeta(\zeta-n)) N_{p}(Q) \label{eq:normalop},
\end{equation}
see the proof of Proposition 5.19 of \cite{mm}, with $\Delta_{p}$ the 
Laplacian on $X_p$ which is the model 
hyperbolic half space up to a linear scaling.  

Now what we are interested in this section, and this paper in general,
is the structure of the difference of the resolvents associated to two
pieces of data.  We begin by proving
\begin{prop}\label{AR.1} Suppose that $g_1, g_2$ are asymptotically 
hyperbolic 
metrics which agree to order $k$ at $\p X,$ i.e in some product decomposition
$X \sim \p X \times [0,\eps)$ near $\p X,$ $x$ is a defining function of 
$\p X,$ in which
\begin{gather}
\begin{gathered}
g_l= \frac{ (dx)^2 +h_l(x,y,dy)}{x^2}
 \\ \text{ where }
h_2(x,y,dy)=h_1(x,y,dy) + x^k L(x,y,dy) + O\left(x^{k+1}\right). \\
\end{gathered}\label{AR.3}
\end{gather}
Suppose  that $V_1, V_2$ are short range
potentials that satisfy $V_1-V_2=x^k W,$ $W \in C^{\infty}(X).$
 Let $\laphd_l,$ $l=1,2,$
be the Laplacian
associated to $g_l$ acting on half-densities via the natural trivialization 
of the half-density bundle given by $g_l$  and let
\begin{gather}
P_l (\zeta) = \laphd_l + V_l + \zeta(\zeta -n).\label{ARPL}
\end{gather}
  Let $h_l(x,y)$ and 
$L(x,y)$  denote the 
matrices of coefficients of the tensors $h_l(x,y,dy)$ and $L(x,y,dy)$
respectively.
We then have that for 
$H=h_1(0,y)^{-1}L(0,y)h_1(0,y)^{-1}$ and 
$T=\operatorname{Tr}\left(h_1(0,y)^{-1}L(0,y)\right)$
\begin{gather}
\begin{gathered}
P_1 - P_2 = x^{k}\left( \sum_{i,j=1}^n H_{ij}x\p_{y_i}x\p_{y_j} 
- \frac{1}{4}k(n-k)T+W\right) +
x^{k+1}R, \\
\end{gathered}\label{AR.2}
\end{gather}
with $R$ a second order symmetric zero-differential operator.
\end{prop}
\begin{pf}
 In local coordinates $(x,y)$ near $q \in \p X$
the operator $P_l$ acts on a $\ha$-density $f(x,y)\ |dx dy|^{\ha}$ as
$$P_l\left(f(x,y) |dx dy|^{\ha}\right)= 
\left[\delta_{l}^{\oq}\left(\Delta_l+V+\zeta(\zeta-n)\right)\delta_{l}^{-\oq} 
f(x,y))\right] |dx dy|^{\ha},$$ where $\Delta_l$ denotes the Laplacian 
acting on functions and $\delta_l$ denotes the determinant of $g_l.$

So we need to consider the operator 
$\delta^{\oq}\Delta_{g}\delta^{-\oq} +V +\zeta(\zeta-n).$ Let
$g_{ij}$ denote the components of the metric $g$  and $g^{ij}$ its 
inverse.  So using the expression of
$\Delta_{g}$ in local coordinates, and denoting $z=(x,y)$ with $z_0=x,$
$z_j=y_j,$ $1 \leq j \leq n,$ we obtain
\begin{gather}
\delta^{\oq}\Delta_{g}\delta^{-\oq}= 
\sum_{i,j=0}^{n} 
\left[-\del^{-\oq}\p_{z_i}\left(g^{ij}\p_{z_j}\del^{\oq}\right) + 
 \del^{-\oq}\p_{z_i}\left( g^{ij}\del^{\oq}\p_{z_j}\right) \right].
 \label{AR.P}
\end{gather}
Recall that
\begin{gather}
\begin{gathered}
g_{00}=\frac{1}{x^2}, \;\ g_{i1}=g_{1i}=0, \;\ i\not = 1, \\
g_{ij}=\frac{1}{x^2}h_{ij}, \;\ i,j \not = 1.
\end{gathered}\label{AR21}
\end{gather}
Therefore
\begin{gather}
\begin{gathered}
g^{00}=x^2, \;\ g^{i1}=g^{1i}=0, \;\ i\not = 1, \\
g^{ij}=x^2h^{ij}, \;\ i,j \not = 1.
\end{gathered}\label{AR211}
\end{gather}
Using \eqref{AR.3} we can write 
\begin{gather}
h_2=h_1\left( I + x^k h_{1}^{-1} L + O\left(x^{k+1}\right)\right)
 \label{AR.4}
\end{gather}
and therefore conclude that
\begin{gather}
\begin{gathered}
h_{2}^{-1}= h_{1}^{-1}  + x^k h_{1}^{-1} L h_{1}^{-1}  + O\left(x^{k+1}\right),
\\
\operatorname{det} h_2= \operatorname{det} h_1 \left( 1 +
x^k\operatorname{Tr}\left(h_{1}^{-1}L\right) + O(x^{k+1})\right).
\end{gathered}\label{AR.5}
\end{gather}
We also deduce from \eqref{AR21} that
\begin{gather}
\delta_j=x^{-2(n+1)} \operatorname{det} h_j, \;\ j=1,2.\label{AR.6}
\end{gather}
Hence
\begin{gather}
\begin{gathered}
\delta_2 =\del_1\left( 1 +
x^k\operatorname{Tr}\left(h_{1}^{-1}L\right) + O(x^{k+1})\right), \\
\delta_2^{\pm \oq} =\del_1^{\pm \oq}\left( 1 \pm
x^k\oq\operatorname{Tr}\left(h_{1}^{-1}L\right) + O(x^{k+1})\right).
\end{gathered}\label{AR.8}
\end{gather}

Examining each term of \eqref{AR.P} and using 
\eqref{AR211}, \eqref{AR.5}, \eqref{AR.6} and \eqref{AR.8} we deduce that
 \eqref{AR.2} holds. This ends the proof of the Proposition.
\end{pf}

Let us denote $P_2-P_1=x^kE,$ where $E$ is the operator given by
the right hand side of \eqref{AR.2}.
Now let $R_1(\zeta)$ be the resolvent of $P_1$ which by Theorem 7.1 of
\cite{mm} lies
in $\Psi DO_{0}^{-2,\zeta,\zeta}(X).$ We then have 
$$P_2 R_1 = P_1 R_1 + (P_2 - P_1) R_1 = \Id + x^k E R_1=P_2R_2+x^k E R_1.$$

So to get $R_2$ as a perturbation of $R_1$ we need to solve,
\begin{equation}
P_2 F = x^{k} E R_1.
\end{equation}
We can rewrite this as,
$$ P_2 ({x'}^{k} F_1) = {x'}^{k} s^{k} E R_1,$$
with $s = x/x'.$ As $x'$ commutes with $P_2,$ this becomes,
$$P_2 F_1 = s^{k} E R_1.$$ Now $s^{k} E R_1$ is in $\Psi
DO_{0}^{0,\zeta+k,\zeta-k}(X),$ we look for $F_1 \in \Psi
DO_{0}^{-2,\zeta,\zeta-k}(X).$ To get improvement on the front face,
we use normal operators,\eqref{AR.P} and the fact that $V$ is short range, to
deduce that
\begin{equation}
(\Delta_p+\zeta(\zeta-n)) N_{p}(F_1) = N_{p} (s^{k} E R_1). \label{normeq} 
\end{equation}
This can be solved near the singularity by using the elliptic calculus,
and away from this the right hand side is in ${\cal
A}^{\zeta+k,\zeta-k}.$ Now Proposition 6.19 of \cite{mm} states that
this equation has a meromorphic solution in ${\cal
A}^{\zeta,\zeta-k}.$  
So we can choose $F_1$ meromorphically to satisfy
\eqref{normeq}. 

Putting $F= {x'}^{k} F_1 \in R^{k}\Psi^{-2,\zeta,\zeta}(X),$ we then have 
that
\begin{gather}
P_{2} (R_1 + F) - \Id \in R^{k+1} \Psi_{0}^{0,\zeta,\zeta}(X).\label{g0}
\end{gather}
We can then remove the term at the front face iteratively and
asymptotically summing obtain, 
$$P_{2}(R_1 + F') - \Id \in R^{\infty} \Psi_{0}^{0,\zeta,\zeta}(X),$$
with
$F' \in R^k \Psi DO^{-2,\zeta,\zeta}(X).$ The error term now vanishes
to infinite order at the front face. The diagonal singularity can be
removed by an element of $ R^{\infty} \Psi_{0}^{-2,\zeta,\zeta}(X)$ by
standard symbolic arguments for constructing the parametrix of a
pseudo-differential operator. This leaves an error in the class
$x^{\zeta} {x'}^{\zeta} f, 
f \in C^{\infty}(X \times X,\Gamma_0^{\ha}(X \times X)).$ This can be
removed using the indicial equation by an element of the same space as
in \cite{mm}. So to summarize, we have proven
\begin{thm}\label{dres}
Let $(X, \p X)$ be a smooth manifold with boundary and defining
function $x.$ Suppose 
$$ g_j = \frac{dx^{2} + h_{j}(x,y,dy)}{x^2}$$ and $V_j$
are smooth real-valued functions vanishing at $\p X.$ Let
$R_{j}(\zeta)$ denote the resolvent of 
$ \laphd_j + V_j + \zeta(\zeta-n)$
where $\laphd_j$ is the Laplacian
associated to $g_j$ acting on half-densities. Suppose that
$\zeta$ is not a pole of $R_j(\zeta).$
Suppose $h_{1} - h_{2}$
and $V_1 - V_2$
vanish to order $k$ at $x=0$ then 
\begin{gather}
R_{1}(\zeta) - R_{2}(\zeta) = G_1(\zeta)+G_2(\zeta) +G_3(\zeta),
\;\ G_i \in \Psi_0^{-2,\zeta,\zeta}(X, \Gamma_0^{\ha}(X)), \;\ i=1,2,3,
\label{dres01}
\end{gather}
where $G_3$ has kernel of the form $x^{\zeta}{x'}^{\zeta}\gamma,$
$\gamma \in C^{\infty}\left(X\times X, \Gamma_0^{\ha}(X\times X)\right),$
the lift of the kernel of $G_2$ under $\beta$ vanishes to infinite order at
the front face  of $X \times_0 X$, and 
 the kernel of $G_1$ satisfies
\begin{gather}
\beta^* G_1(\zeta)= R^{k}\rho^{\zeta} {\rho'}^{\zeta} \alpha(\zeta), \;\
\alpha(\zeta) \in 
C^{\infty}\left(X\times_0 X \setminus \Delta_0, 
\Gamma_0^{\ha}(X\times_0 X)\right), \label{dres1}
\end{gather}
is a conormal distribution to the lifted diagonal $\Delta_0.$

If $E$ is such that $Q_1 - Q_2 = x^{k} E,$ then the restriction of 
$\alpha(\zeta)$ to the front face satisfies
\begin{gather}
(\Delta_{h_0} + \zeta(\zeta-n))
N_{p}\left(\rho^{\zeta}{\rho'}^{\zeta-k}\alpha(\zeta) 
\right) 
= N_{p}(\fdl{x}{x'}^{k}E)G, \label{dres02}
\end{gather}
 where $F$ is the front face, $\{R=0\},$ $\Delta_{h_0}$ is the Laplacian on 
the hyperbolic  space with metric $h_{0}(p),$ i.e in coordinates $(z_0,z')$
where the boundary is $\{z_0=0\},$
\begin{gather}
\Delta_{h_0}= z_{0}^{2} \sum_{i,j=0}^{n} h^{ij}(p) \p_{z_i}\p_{z_j}-
(n-1)z_0\p_{z_0}, \label{deltah}
\end{gather}
 and $G$ is the Green's function of $\Delta_{h_0}+\zeta(\zeta-n)$. 
\end{thm}
Note that the last statement follows from Propositions 2.17 and 5.19 of
\cite{mm} and the fact that the normal operator of the resolvent is its 
Green's function.
\begin{remark}\label{remark1}
In what follows it is important to realize that there is a unique
solution 
of \eqref{dres02}  which is meromorphic in $\zeta,$ is conormal to the 
centre of  the front face and such that near the boundaries is in
${\mathcal A}^{\zeta,\zeta-k}.$ 

  To see that note that, if we  
have two choices, $w_1$ and $w_2,$
then $(\Delta_{h_0} + \zeta(\zeta-n))(w_1-w_2)=0.$
Since $w_1-w_2$ is conormal to the centre of the front face it must be 
actually smooth. 
By Theorem 7.3 of \cite{maz1},
we know that $w_1-w_2=(\rho\rho')^{\zeta} f + (\rho\rho')^{n-\zeta} g$
where  $f,g$ are distributional coefficients. 

On the other hand
we also know that  $w_1-w_2 \in {\mathcal A}^{\zeta,\zeta-k},$
 and so is of the form
$w_1-w_2=\rho^{\zeta}{ \rho'}^{\zeta-k} w$ with $w$ smooth up to the boundary. 
Therefore we conclude that 
$w_1-w_2=(\rho \rho')^{\zeta}\widetilde{w}$ with $\widetilde{w}$ smooth up to
the boundary. Since $\Delta_{h_0}$ has no discrete spectrum,
it follows from Proposition \ref{prop:eiguniqueness} that for
$\zeta \not \in (-\infty,\frac{n}{2}],$ $w_1=w_2,$ and thus by
meromorphicity everywhere. 
\end{remark}

\section{The Poisson Operator and The Scattering Matrix}\label{scatmat}

In this section we extend some of the results of \cite{guillop} and 
\cite{gz1}, obtained in the case of Riemann surfaces, to asymptotically 
hyperbolic manifolds.  We show that the kernel of the Poisson operator is
a multiple of the Eisenstein function and, as in \cite{guillop}, we
obtain a formula for the scattering matrix in terms of the resolvent. (Similar
results have been established by Borthwick in \cite{Bo} for $\Re \zeta=n/2.$)
As a consequence of this formula, we prove that
the scattering matrix at energy $\zeta,$ $\zeta \in \mc \setminus Q,$
where $Q$ is a discrete subset which is described in Proposition \ref{PO7},
 is a pseudo-differential operator of order $2\zeta-n.$
We also prove the result stated in equation \eqref{sca1} 
of the introduction.  

Before proceeding to this, we sketch our argument. The resolvent of
the Laplacian acting on half-densities has by Mazzeo-Melrose, \cite{mm}, a
meromorphic extension to the entire complex plane. Its weighted restriction to 
$X \times \p X,$ we call the Eisenstein function, $E(\zeta),$ in analogy to
previous work on hyperbolic manifolds. This function is automatically
in the kernel of $\Delta+\zeta(\zeta-n)$ and we examine its
distributional asymptotics. In particular, we see that they have two
components one lead term is a multiple of the 
delta function on the diagonal times $x^{\zeta}$ and the other is a 
pseudo-differential operator times $x^{2n-\zeta}.$ This means that
upon integration of a suitable multiple of 
the Eisenstein function against a half-density on the boundary one
obtains roughly 
an eigenfunction of the form $x^{\zeta} f + x^{2n-\zeta}g$ plus lower
order terms, 
where $g=S(\zeta)f$ with $f$ prescribed and $S(\zeta)$ a fixed
pseudo-differential 
operator which is of course the scattering matrix acting on
half-densities. So the Eisenstein function is really the Poisson
operator for the problem and our first task is to prove it has the
appropriate distributional asymptotics.  The Eisenstein function
$E(\zeta)$ plays an analogous role to that of the Poisson operator
$P(\lambda)$ in \cite{smagai}. However it lives on the manifold, $X
\times \p X$
blown-up along the boundary diagonal rather than on a micro-locally
blown-up space. 

Recall that $X\times_0 X$ is the space obtained from $X\times X$ by 
blowing-up
the diagonal $\Delta \subset \p X \times \p X$ and that
$\beta: X\times_0 X \longrightarrow X\times X$
is the corresponding blow-down map.
Theorem 7.1 of \cite{mm}  states that
the resolvent $R(\zeta)=\left( \laphd_g+V+\zeta(\zeta-n)\right)^{-1},$
which is well defined for $\Re\zeta$ large, extends to a meromorphic family
$R(\zeta)\in \Psi_0^{-2,\zeta,\zeta}(X,\Gamma_0^{\ha}(X))$ 
$\zeta \in \mc,$ that satisfies,
in terms of the spaces introduced in section \ref{S.res},
\begin{gather*}
R(\zeta)=R'(\zeta)+R''(\zeta), \;\ 
R'(\zeta) \in \Psi_0^{-2}(X,\Gamma_0^{\ha}(X)) \text{ and }
R''(\zeta) \in \Psi_0^{-\infty,\zeta,\zeta}(X,\Gamma_0^{\ha}(X)),
\end{gather*}
with the boundary term, $R''(\zeta),$ having Schwartz  kernel of a special 
form
\begin{gather}
\beta^* K''(\zeta)=\rho^{\zeta}{\rho'}^{\zeta}F(\zeta), \;\
F(\zeta) \in C^\infty\left(X \times_0 X; 
\Gamma_0^{\ha}(X\times_0 X)\right)\label{PO1}
\end{gather}
where $\rho$ and $\rho'$ are defining functions of the 
top and bottom faces respectively, and $F(\zeta)$ is meromorphic in $\zeta.$

Let $R(\zeta) \in C^{-\infty}(X \times X, \Gamma_0^{\ha}(X\times X))$ 
also denote the Schwartz kernel of the resolvent and let
 $x$ and $x'$ be boundary defining function of each copy of $X$ in 
$X\times X.$   We will show that the Eisenstein function, which is
defined by
\begin{gather}
E(\zeta)= {x'}^{-\zeta+\frac{n}{2}} R(\zeta)|_{x'=0}, \label{PO2}
\end{gather} 
is a smooth section of $\Gamma_0^{\ha}(X \times \p X).$ 
Notice that it 
depends on the choice of the defining function $x'.$ To make it independent of
this choice one can view it as a section of
$\Gamma_0^{\ha}(X \times \p X)\otimes |N^*\p X|^{\zeta-\frac{n}{2}}$
by defining it as
\begin{gather}
E(\zeta)= {x'}^{-\zeta+\frac{n}{2}} R(\zeta)|_{x'=0}
|dx'|^{\zeta-\frac{n}{2}}. \label{PO21}
\end{gather} 
This is the analogue of Definition 2.2 of \cite{guillop}.  For simplicity
we will work with the definition given by \eqref{PO2} and so  we fix a product
decomposition $X \sim \p X \times [0,\eps)$ of $X$ near $\p X.$

Since $R'(\zeta) \in \Psi_0^{-2}(X,\Gamma_0^{\ha}(X)),$ its kernel vanishes 
to infinite order at the top and bottom faces.  So
we deduce that its kernel satisfies 
\begin{gather*}
{x'}^{-\zeta+\frac{n}{2}} R'(\zeta)|_{x'=0}=0.
\end{gather*}
Therefore, for $K''$ given by \eqref{PO1},
\begin{gather}
E(\zeta)= {x'}^{-\zeta+\frac{n}{2}}K''(\zeta){|_{x'=0}}. \label{PO3}
\end{gather} 

Next we blow-up the manifold $X \times \p X$ along $\Delta$ and analyze the
lift of $E(\zeta)$ under the blow-down map.
Let $X \times_{0} \p X$ be the manifold with corners obtained by 
blowing-up $X \times \p X$ along the diagonal
$\Delta \subset \p X \times \p X$ and let
\begin{gather*}
\widetilde{\beta} : X \times_{0} \p X \longrightarrow X \times \p X
\end{gather*}
denote the corresponding blow-down map.  It is then clear that
 $\widetilde{\beta}= \beta|_{(X\times_0 \p X)}.$

 Let $\mcf$ be the new boundary
face introduced by the blow-up, the front face, i.e
\begin{gather*}
{\mcf}={\widetilde{\beta}}^{-1}\left(\Delta\right).
\end{gather*}
If $R \in C^\infty(X \times_0 X)$ is a defining function of the front face
in $X\times_0 X,$
$R_{|_{X \times_0 \p X}}$ is a defining function of $\mcf,$  that we will
also denote by $R.$

Next we consider the lift of $E$ under the map $\widetilde{\beta}.$ It is
actually more convenient to analyze the lift of 
$x^{-\zeta+\frac{n}{2}}E$ first.
We deduce from \eqref{PO1} and \eqref{PO3} that

\begin{gather}
{\widetilde{\beta}}^*\left(x^{-\zeta+\frac{n}{2}}E(\zeta)\right)= 
(R\rho)^{-2\zeta+\frac{n}{2}+\ha}
(R\rho')^{-\zeta+\frac{n}{2}}
\rho^{\zeta} {\rho'}^{\zeta}F(\zeta){|_{\rho'=0}}=
R^{-2\zeta +\frac{n}{2}} \rho^{\frac{n}{2}}
{\rho'}^{\frac{n}{2}}F(\zeta){|_{\rho'=0}}. \label{PO3O} 
\end{gather}

As in section 3 of \cite{mm}, in the  region 
away from the bottom face,  we can use projective coordinates
$x,$ $\rho'=x'/x,$ $Y=(y-y')/x.$  So  we can represent the half-density 
$F(\zeta)$ in these local coordinates by
\begin{gather}
 F(\zeta)= F(\zeta,\rho',x,y,Y) \left|\frac{d x}{x}\frac{dy}{x^n}
\frac{d \rho'}{\rho'}\frac{d Y}{{\rho'}^n}\right|^{\ha}.\label{PO4}
\end{gather}

Near the intersection of the top and bottom faces we can use local
coordinates $R=|y-y'|,$ $\rho=x/R,$ $\rho'=x'/R$ and $Y=(y-y')/R.$
Then $F(\zeta)$ can be represented by
\begin{gather}
 F(\zeta)= F(\zeta,\rho,\rho',R,y,Y) \left|\frac{d R}{R}\frac{d \rho'}{\rho'}
\frac{d \rho}{\rho}\frac{d Y dy}{\rho^n{\rho'}^nR^n}\right|^{\ha}.
\label{PO41} 
\end{gather}

Thus it follows from \eqref{PO3O} that
${\widetilde{\beta}}^*(x^{-\zeta+\frac{n}{2}+\ha}E)$ is given in these 
local coordinates respectively by
\begin{gather}
x^{-2\zeta+\frac{n-1}{2}}F(\zeta,0,x,y,Y) 
\left|dx dy dY\right|^{\ha}, \;\ 
R^{-2\zeta+\frac{n-1}{2}}\rho^{-\ha}F(\zeta,\rho,0,R,y,Y) 
\left|d R d \rho d Y dy \right|^{\ha}.\label{4.7new}
\end{gather}
Therefore we have that
\begin{gather}
{\widetilde{\beta}}^*(x^{-\zeta+\frac{n}{2}}E) \in 
R^{-2\zeta+\frac{n-1}{2}}\rho^{-\ha}
C^\infty(X \times_0 \p X, \Gamma^{\ha}(X \times_0 \p X)), \label{POgb}
\end{gather}

 Notice that 
\begin{gather}
\widetilde{\beta}^* :x^{-\ha}\Gamma^{\ha}(X \times \p X)
\longleftrightarrow R^{\frac{n}{2}}\rho^{-\ha}\Gamma^{\ha}(X \times_0 \p X) 
\label{pb1}
\end{gather}
is an isomorphism.  To see that we use
local coordinates 
\begin{gather}
\begin{gathered}
x,\;\  Y=(y-y')/x, \;\ y \\
|y-y'|=R, \;\ \rho=x/R, \;\ Y=(y-y')/R
\end{gathered}\label{PFloc}
\end{gather}
where the first set is valid away from 
$M=\operatorname{clos}{\widetilde{\beta}}^{-1}
\left(\p X \times \p X \setminus \Delta\right)$  and the second is valid
near $M \cap \mcf$ respectively. Then the lift of
$x^{-\ha}\left|dx dy dy'\right|^{\ha}$ is given by
\begin{gather*}
x^{\frac{n-1}{2}} \left|dx dy dY\right|^{\ha}, \;\ \;\
R^{\frac{n-1}{2}}\rho^{-\ha} \left|dR d\rho dY dy\right|^{\ha}
\end{gather*}
respectively.
Therefore the map \eqref{pb1} is in fact an isomorphism.

Next we consider the push-forward of a smooth section  of 
$R^{-2\zeta+\frac{n-1}{2}}\rho^{-\ha} \Gamma^{\ha}(X \times_0 \p X).$ 
First we need to introduce some notation.
Note that $\mcf,$ $M$ are  manifolds with boundary,
 and that the restriction of $\widetilde{\beta}$ to $M$ induces
a map 
\begin{gather*}
\beta_{\p}=\beta|_{M} : M \sim \p X \times_0 \p X \longrightarrow 
\p X \times \p X
\end{gather*}
which corresponds to the blow-up of the manifold $\p X \times \p X$ 
along the diagonal $\Delta \subset \p X \times \p X.$

Given  $R \in C^\infty( X \times_0 \p X)$ and $x \in C^\infty(X),$ defining
functions of ${\mcf}$ and $\p X$ respectively, the function
$\rho=\frac{x}{R} \in C^\infty( X \times_0 \p X)$
is a defining function of $M.$  Since ${\mcf}$ and $M$ intersect
transversally, with $M\cap \mcf=\p M=\p \mcf,$ the functions
\begin{gather*}
R_M= R|_M \in C^\infty(\p X \times \p X), \;\ \;\ 
\rho_{\mcf}= \rho|_{\mcf} \in C^\infty(\mcf)
\end{gather*}
are defining functions of  $\p M$  and $\p {\mcf}$ 
respectively.

Recall that, see for example section 3.2 of \cite{Ho1}, if $Y$ is a 
manifold
with corners and $y \in C^\infty(Y),$ is a defining function of a 
boundary hypersurface of $Y,$ then  sections of 
$y^{\zeta} \Gamma^{\ha}(Y),$ viewed as distributions acting on
$\Gamma^{\ha}(Y)$ via
\begin{gather}
\begin{gathered}
\lan y^{\zeta}F, f\ran = \int_Y y^\zeta Ff, \;\ 
\text{ for } \Re \zeta >-1, \;\ 
 F \in C^\infty(Y, \Gamma^{\ha}(Y)), \;\
 f \in C_0^\infty(Y, \Gamma^{\ha}(Y)), 
\end{gathered}\label{PF2}
\end{gather}
have holomorphic extensions to $\mc\setminus -\mn.$ 

We will consider three such half-densities associated to
$R,$ $R_M$ and $\rho_{\mcf}$ defined on $X\times_0 \p X,$
$M,$ and $\mcf$ respectively.

We have fixed a product decomposition $X \sim \p X \times [0,\eps),$
near $\p X,$ and  will prove that the sections of the push-forward of
$R^{-2\zeta+\frac{n-1}{2}}\rho^{-\ha} \Gamma^{\ha}(X \times_0 \p X)$  have
distributional asymptotic expansions as $x \downarrow 0.$ 
To do that we define the partial pairing for
$u \in
 R^{-2\zeta+\frac{n-1}{2}}\rho^{-\ha}C^{\infty}(X\times_0 \p X,
 \Gamma^{\ha}(X \times_0 \p X)),$
$f \in C^{\infty}(\p X \times \p X, \Gamma^{\ha}(\p X \times \p X))$

\begin{gather}
\lan {\beta_{\p}}_*u, f\ran = \int_{\p X \times \p X} 
\left({\beta_{\p}}_* u\right)(x,y,y') f(y,y').\label{ppar}
\end{gather}

We remark that if $u$ is a smooth section of
$R^{-2\zeta+\frac{n-1}{2}}\rho^{-\ha} \Gamma^{\ha}(X \times_0 \p X)$ then
the restriction of $u$ to $M,$ denoted by $u|_{M}$ is well defined as
a section of 
$R^{-2\zeta+\frac{n-1}{2}}\rho_{\mcf}^{-\ha} \Gamma^{\ha}(\p X \times_0 \p X).$
It is also easy to see that 
$x^{2\zeta-\frac{n-1}{2}} u$ is a smooth section of 
$\rho^{2\zeta-\frac{n}{2}} \Gamma^{\ha}(\p X \times_0 \p X).$ Therefore
it can be restricted to $\mcf=\{R=0\}$ and 
$(x^{2\zeta-\frac{n}{2}} u)|_{\mcf}$ is a smooth section of
$\rho_{\mcf}^{2\zeta-\frac{n}{2}} \Gamma^{\ha}(\mcf).$

We now prove a push-forward theorem which relates the distributional
asymptotics of a class of half-densities including the Eisenstein function
to their behaviour at the boundary, cf Prop 16 of \cite{smagai}.
\begin{prop}\label{PF3} Let $x\in C^\infty(X)$ be a defining function of 
$\p X$  and fix  a product decomposition $X \sim \p X \times [0,\eps)$ 
near $\p X.$ Let 
$R \in C^\infty(X \times_0 \p X)$ be a defining 
function of $\mcf,$ 
and let $\rho_{\mcf}$  be defined as above. Let 
$v=R^{-2\zeta+\frac{n-1}{2}} \rho^{-\ha}F,$
 $F \in C^{\infty}(X \times_0 \p X, \Gamma^{\ha}(X \times_0 \p X)),$ 
$2\zeta \in \mc \setminus \mz.$ 
Then the push-forward of $v$ under $\widetilde{\beta},$  denoted by
$\widetilde{\beta}_{*}v,$ is a section of 
$x^{-\ha}\Gamma^{\ha}(X \times \p X),$ which
has  a conormal singularity at $\Delta,$ and moreover it has an asymptotic 
expansion in $x$ as $x\downarrow 0,$ in the sense that if 
$f \in C^{\infty}(\p X \times \p  X, \Gamma^{\ha}(\p X \times \p  X)),$ and 
$\lan\; ,\; \ran$ is the partial
pairing defined above, then 
\begin{gather}
\lan\widetilde{\beta}_{*}v, f\ran =(H_{\zeta}(x)+ 
x^{n-2\zeta}G_{\zeta}(x))\fracwithdelims||{dx}{x}^{\ha} , \text{ as } x\downarrow 0, \label{PF4}
\end{gather}
where $G_\zeta, H_\zeta \in
C^{\infty}([0,\eps))$ 
depend holomorphically on $\zeta.$ Moreover if $v|_M$ and 
$x^{2\zeta-\frac{n-1}{2}}v|_{\mcf}$ denote the restrictions of these 
half-densities to $M$ and $\mcf$ respectively, then
\begin{gather}
\begin{gathered}
H_\zeta(0)= 
\lan {\beta_{\p}}_*\left(v|_{M}\right), f\ran, \\
G_\zeta(0)= \left\lan \lan x^{2\zeta-\frac{n-1}{2}}v|_{\mcf},\rho_{\mcf}^{-\frac{n+2}{2}} \ran \delta_\Delta, f\right\ran 
\end{gathered}\label{PF5}
\end{gather}
where $\delta_{\Delta}$ is the delta function of the diagonal and
$\lan x^{2\zeta-\frac{n-1}{2}}v|_{\mcf},\rho_{\mcf}^{-1-\frac{n}{2}}\ran$
is the pairing induced by the trivialization of the half-density
bundle $\Gamma^{\ha}({\mcf})$ given by the product structure.
\end{prop}
\begin{pf}  
Since this is a local result and $\widetilde{\beta}$ is a
diffeomorphism away from $\Delta,$ we only need to work in a neighbourhood of
a point $q \in \Delta.$  Let $y, y'$ be local coordinates near $q$ and let
$R=\left(x^2+|y-y'|^2\right)^{\frac{1}{2}},$ $\rho=x/R$ and
$\omega=(y-y')/R.$ The map $\widetilde{\beta}$ can be described as
\begin{gather*}
\widetilde{\beta}: {\ms}_{+}^{n} \times [0,\infty) \times \mrn 
\longrightarrow \mr \times {\mr}^{n-1}  \times {\mr}^{n-1} \\
(\rho,\omega,R) \longmapsto (R\rho, y'+R\omega, y').
\end{gather*}
and we will denote
\begin{gather}
v= R^{-2\zeta+\frac{n-1}{2}}\rho^{-\ha}F(\rho, \omega, R)
\left|d\rho d\omega dR\right|^{\ha}, \; F \in C^{\infty}(X \times_0 \p X).
\label{PFv} 
\end{gather}

We also set $z=y-y'.$ Then the variables $y'$ become parametric and 
for simplicity we will ignore them.  The diagonal is given by
$\Delta=\{x=0, z=0\}.$ First we observe that
the vector fields tangent to $\Delta$ are spanned over 
$C^\infty(X\times \p X)$ by
\begin{gather*}
x\p_x, \;\ z_j\p_{z_k}, \;\ x\p_{z_k}, \;\ z_k\p_x
\end{gather*}
and it can be proven, by using projective coordinates as in \eqref{PFloc} 
above,
that these vector fields lift  under $\widetilde{\beta}$
to smooth vector fields that are tangent to $\mcf.$ 
This shows that $\widetilde{\beta}_{*}(R^{-2\zeta}F)$ is conormal to 
$\Delta.$

We observe that the radial vector field is given by
\begin{gather}
R\p_R={\widetilde{\beta}}^*\left(x\p_x + z\cdot\p_z \right).
\label{PF6}
\end{gather}
Thus, since $x(x^2+|z|^2)^{-\frac{1}{2}}$ is homogeneous of degree zero with 
respect to the action $(x,z) \mapsto (\la x, \la z),$ $\la \in \mr_+,$ 
we have
\begin{gather*}
(x\p_x + z\cdot\p_z)\left(x(x^2+|z|^2)^{-\frac{1}{2}}\right)=0, \;\ 
(x,z)\not=0.
\end{gather*}
Therefore, 
\begin{gather}
R\p_R\frac{x^k}{R^k}
 =0, \;\ k \in \mn. \label{PF7} 
\end{gather}

We will also use that

\begin{gather}
\p_x (x^2 + |z|^2)^{\frac{1}{2}} = x(x^2 + |z|^2)^{-\frac{1}{2}}
\label{PF8}
\end{gather}
and that
$R{\widetilde{\beta}}^*(\p_x)$ is a smooth vector field in $X\times_0 \p X.$

We deduce from \eqref{PF6} and \eqref{PF8} that
\begin{gather*}
x\p_x (x\p_x + z\cdot \p_z + 2\zeta){\widetilde{\beta}}_*(R^{-2\zeta}F)=
x{\widetilde{\beta}}_*(R^{-2\zeta}F_1), \;\ F_1 \in 
C^{\infty}(X\times_0 \p X)
\end{gather*}
and using \eqref{PF6}, \eqref{PF7} and \eqref{PF8}  we obtain
\begin{gather*}
(x\p_x-1) (x\p_x + z\cdot \p_z + 2\zeta-1)
x{\widetilde{\beta}}_*(R^{-2\zeta}F_1)=
x^2{\widetilde{\beta}}_*(R^{-2\zeta}F_2),
 \;\ F_2 \in C^{\infty}(X\times_0 \p X).
\end{gather*}

Similarly we obtain
\begin{gather}
\begin{gathered}
(x\p_x-k) (x\p_x + z\cdot \p_z + 2\zeta-k)
x^k{\widetilde{\beta}}_*(R^{-2\zeta}F_k)=
x^{k+1}{\widetilde{\beta}}_*(R^{-2\zeta}F_{k+1}), \\
F_{k+1} \in C^{\infty}(X\times_0 \p X).
\end{gathered}\label{PF81}
\end{gather}
By induction and  \eqref{PF81} we obtain

\begin{gather}
\begin{gathered}
\prod_{j=0}^{M}(x\p_x-j)(x\p_x+z\cdot\p_z+2\zeta-j)
{\widetilde{\beta}}_*(R^{-2\zeta}F)=
x^{M+1}{\widetilde{\beta}}_*(R^{-2\zeta}F_M), \\
F_M \in C^{\infty}(X\times_0 \p X). 
\end{gathered}\label{PF11}
\end{gather}

Since the map defined in \eqref{pb1} is an isomorphism, it follows that
the push-forward of \eqref{PFv}
can be written in local coordinates $x,z$ as
\begin{gather*}
{\widetilde{\beta}}_*(R^{-2\zeta}F)(x,z)\left|\frac{dx}{x} dz\right|^{\ha}.
\end{gather*}
Let $f \in C^\infty(\mrn).$ Then 
\begin{gather*}
\lan {\widetilde{\beta}}_*(v), f(z)\left|dz\right|^{\ha}\ran = 
\left(\int_{\mrn} {\widetilde{\beta}}_*(R^{-2\zeta}F)(x,z)f(z) dz \right)
\left|\frac{dx}{x} \right|^{\ha}.
\end{gather*}

Using \eqref{PF11} and the identity 
$\operatorname{div}(zu(z))= n u(z)+z\cdot\p_zu(z)$ we deduce that 
\begin{gather*}
u(x)=\int_{\mrn} {\widetilde{\beta}}_*(R^{-2\zeta}F)(x,z)f(z) dz,
\;\ x>0
\end{gather*}
satisfies
\begin{gather}
\begin{gathered}
\prod_{j=0}^{M} \left(x\p_x-j\right)\left(x\p_x+2\zeta-(n+j)\right)u(x)=
x^{M+1}\int_{\mrn} {\widetilde{\beta}}_*(R^{-2\zeta}F_M)(x,z)f(z)dz, \;\
x>0.
\end{gathered} \label{PF12}
\end{gather}
Let $$u_M(x)=\prod_{j=0}^{M}\left(x\p_x+2\zeta-(n+j)\right)u(x).$$ Then we 
deduce from \eqref{PF12} that there exists $s \in \mr,$ independent of $M,$
such that
\begin{gather}
\left|\p_x\prod_{j=1}^{M}\left(x\p_x-j\right)u_M\right|
 \leq C|x|^{M+s}, \;\ x>0 .\label{PFB1}
\end{gather}
Thus, for $M$ large, there exists $C_0 \in \mc$ such that
\begin{gather}
\lim_{x\downarrow 0} \prod_{j=1}^{M}\left(x\p_x-j\right)u_M=C_0. \label{PFB2}
\end{gather}
From \eqref{PFB1} we obtain
\begin{gather*}
\left|\p_x\left(\prod_{j=1}^{M}\left(x\p_x-j\right)u_M-C_0\right)\right|
 \leq C|x|^{M+s}, \;\ x>0.
\end{gather*}
It follows from \eqref{PFB2} that
\begin{gather*}
\left|\prod_{j=1}^{M}\left(x\p_x-j\right)u_M-C_0\right| \leq C|x|^{M+s+1},
\;\ x>0.
\end{gather*}
Since $(x\p_x-j)C_0=-jC_0$ we obtain, for $a_0=C_0/(-1)^M M! ,$
\begin{gather*}
\left|\prod_{j=1}^{M}\left(x\p_x-j\right)(u_M-a_0)\right| \leq C|x|^{M+s},
\;\ x>0.
\end{gather*}
Proceeding by induction we find that for $M+s-p+1>0,$ there exist 
$a_m\in \mc,$ $0\leq m \leq p-1,$ depending on $\zeta,$ such that
\begin{gather}
\left|\prod_{j=p}^{M}\left(x\p_x-j\right)
(u_M-\sum_{m=0}^{p-1}a_m x^m)\right| \leq C|x|^{M+s+1}, \;\ x>0.\label{PFA5}
\end{gather}
Now we observe that if
\begin{gather*}
|(x\p_x-\alpha)u| \leq C|x|^\beta, \;\ x>0, \;\ \beta < \Re \alpha
\end{gather*}
then 
\begin{gather}
|u(x)| \leq C|x|^\beta.\label{PFA6}
\end{gather}
Indeed, just notice that
\begin{gather*}
|u(1)-x^{-\alpha}u(x)|=\left|\int_{x}^1 \p_s(s^{-\alpha}u(s))ds\right|=
\left|\int_{x}^1 s^{-\alpha-1}\left(s\p_s-\alpha\right)u(s) ds\right| \leq \\
C\int_{x}^1 s^{-\alpha-\Re \beta -1}ds=
\frac{C}{\beta-\Re \alpha}\left(1-x^{\beta-\Re \alpha}\right).
\end{gather*}
Thus \eqref{PFA6} follows.  Therefore we deduce from \eqref{PFA5} and
\eqref{PFA6} that,
for $M+s-p>0,$
\begin{gather}
\left| u_M(x)-\sum_{m=0}^{p-1} a_mx^m\right| \leq C|x|^{M+s}, \;\
x>0.\label{PFA7}
\end{gather}
Notice that $(x\p_x+2\zeta-n-j)x^m=(m+2\zeta-n-j)x^m.$ Since 
$2\zeta\not\in \mz,$ $j\not=2\zeta-n-j,$ so
 we deduce from \eqref{PFA7} that for 
$d_m\prod_{j=0}^M (m+2\zeta-n-j)=a_m,$  $v_p=\sum_{m=0}^{p-1} d_mx^m$
satisfies
\begin{gather}
\left| \prod_{j=0}^M \left(x\p_x +2\zeta-n-j\right)\left(u-v_p\right) \right|
\leq C|x|^{M+s+1}. \label{PFA8}
\end{gather}
This gives that
\begin{gather*}
\left|x^{-2\zeta+n+1}\p_x x^{2\zeta-n} \prod_{j=1}^M 
\left(x\p_x +2\zeta-n-j\right)\left(u-v_p\right) \right| \leq C|x|^{M+s+1}. 
\end{gather*}
Thus, for $M+s+2\Re \zeta-n-2>0,$ there exists $b_0 \in \mc$ such that
\begin{gather*}
\lim_{x\downarrow 0} x^{2\zeta-n} \prod_{j=1}^M 
\left(x\p_x +2\zeta-n-j\right)\left(u-v_p\right)=b_0.
\end{gather*}
Since $2\zeta\not \in \mz,$ we can proceed as above to deduce that
there exists $\gamma_0$ such that
\begin{gather*}
\left| \prod_{j=1}^M
 \left(x\p_x +2\zeta-n-j\right)\left(u-v_p-\gamma_0 x^{n-2\zeta}\right)\right|
\leq C|x|^{M+s+1}, \;\ x>0.
\end{gather*}
Using induction we find that for $M+s+2\Re \zeta-n-q>0,$
there exist $\gamma_m,$ $0\leq m\leq q-1,$ depending on $\zeta,$ such that
\begin{gather*}
\left| \prod_{j=q}^M
 \left(x\p_x +2\zeta-n-j\right)\left(u-v_p-\sum_{m=0}^{q-1} 
\gamma_m x^{m+n-2\zeta}\right)\right| \leq C|x|^{M+s+1}, \;\ x>0.
\end{gather*}
From \eqref{PFA6} we obtain, for arbitrary $M\in \mn,$
and $p,q$ satisfying respectively $M +s -p>0,$ $M+s+2\Re \zeta-n-q>0,$
\begin{gather*}
\left|u(x)- \sum_{m=0}^{p-1}d_m x^m -
x^{n-2\zeta}\sum_{m=0}^{q-1}\gamma_m x^m \right|
\leq C|x|^{M+s+1}, \;\ x>0
\end{gather*}
Now Borel's lemma gives the desired result. 
It is clear from the construction 
that $H_{\zeta}$ and $G_\zeta$ depend holomorphically in $\zeta,$ provided
$2\zeta \not \in \mz.$  This method of proving the existence of
an expansion goes back to Euler and has been used in similar contexts
in \cite{radial}, \cite{forward} and also \cite{maps}.

Next we need to compute $G_\zeta(0)=\gamma_0$ and $H_\zeta(0)=d_0.$ Since 
these are
holomorphic functions of $\zeta,$  we only need 
to compute $H_\zeta(0)$ for $2\Re\zeta-n>0,$ and $G_\zeta(0)$ for
$2\Re\zeta-n<0$. 

In the coordinates above we have, for $f\in C_0^\infty(\mrn)$ and
$R=(x^2+|z|^2)^{\frac{1}{2}},$
\begin{gather}
\lan{\widetilde{\beta}}_{*}v, f|dz|^{\ha} \ran = 
\int_{\mrn} R^{-2\zeta} F\left(R,
\frac{x}{R}, \frac{z}{R}\right) f(z) dz |dx|^{\ha}\label{PFA9}
\end{gather}
It follows from the dominated convergence theorem that for $2\Re\zeta-n>0,$
\begin{gather*}
\lim_{x\downarrow 0}  \int_{\mrn} R^{-2\zeta}F\left(R,
\frac{x}{R}, \frac{z}{R}\right) f(z) dz= \int_{\mrn} |z|^{-2\zeta}
 F\left(|z|,0,\frac{z}{|z|}\right)f(z) dz= \\
\lan{\beta_{\p}}_* \left(v |_M\right),f\ran.
\end{gather*}
To compute $G_\zeta(0)$ for $2\Re\zeta-n <0$ we set $z=xw.$  
Observing that in these coordinates $\rho_{\mcf}=(1+|w|^2)^{-\frac{1}{2}},$
we deduce from \eqref{PFA9} that 
\begin{gather*} 
\lan{\widetilde{\beta}}_{*}(R^{-2\zeta} F),f\ran = 
x^{n-2\zeta} \int_{\mrn} \rho_{\mcf}^{2\zeta}F\left(x\rho_{\mcf}^{-1},
 \rho_{\mcf}, \frac{w}{(1+|w|^2)^{\frac{1}{2}}}\right) f(xw) dw.
\end{gather*}
Again by the dominated convergence theorem
\begin{gather*}
\lim_{x\downarrow 0}  
\int_{\mrn} \rho_{\mcf}^{2\zeta}F\left(x\rho_{\mcf}^{-1},
 \rho_{\mcf}, \frac{w}{(1+|w|^2)^{\frac{1}{2}}}\right) f(xw) dw= 
f(0) \int_{\mrn} \rho_{\mcf}^{2\zeta}
F\left(0, \rho_{\mcf}, \frac{w}{(1+|w|^2)^{\frac{1}{2}}}\right) dw.
\end{gather*}

Using the map $w \ni \mrn \longmapsto 
SP(w)=\left( (1+|w|^2)^{-\ha},w(1+|w|^2)^{-\ha}\right) \in \msn$ we have,
for $g \in C^\infty(\msn_{+}),$ 
$$\int_{\ms_{+}^{n}}g d\sigma= \int_{\mrn} g(SP(w))(1+|w|^2)^{-\frac{n+1}{2}} 
dw.$$
Therefore
\begin{gather}
\begin{gathered}
\int_{\mrn} \rho_{\mcf}^{2\zeta}
F\left(0, \rho_{\mcf}, \rho_{\mcf}w\right) dw=
\int_{\msn_{+}}\rho_{{\mcf}}^{-1-\frac{n}{2}}\rho_{\mcf}^{2\zeta-\frac{n}{2}} 
F d\sigma= \\
\lan x^{2\zeta-\frac{n-1}{2}}v|_{\mcf},\rho_{\mcf}^{-1-\frac{n}{2}}\ran.
\end{gathered}\label{1018}
\end{gather} 
This concludes the proof of the Proposition.
\end{pf}

The following Proposition will be important in the definition of the
scattering matrix.
\begin{prop}\label{G01} Let $g$ and $V$ satisfy the hypotheses of Theorem 
\ref{pssm0}.  The coefficient of $\delta_{\Delta}$ in the second equation of
\eqref{PF5}, 
$M(\zeta)=
\lan
x^{2\zeta-\frac{n-1}{2}}E(\zeta)|_{\mcf},\rho_{\mcf}^{-1-\frac{n}{2}}\ran,$
is equal to $|h_{0}|^{\ha}$ times a function which 
is independent of the base point of the fibre $\mcf$ and is also
 independent of $g$ and $V.$
\end{prop}
\begin{pf}
According to \eqref{1018}, \eqref{PO3O} and \eqref{pb1}, $M(\zeta)$ depends 
only on the
value of $F|_{\mcf}$ where $F$ is given by \eqref{PO1} and $\mcf$ is
as above.

We recall from the construction of $R(\zeta)$ in section 3 and the
proof of Proposition 7.4 in \cite{mm} that the normal operator of the
$R(\zeta)$ is just $R_{0}(\zeta),$ the Green's function of the operator
$\Delta_{h_0}+\zeta(\zeta-n)$ given by \eqref{deltah}, where as 
observed  in \cite{mm}, the fibre of the
front face over a point $p \in \p M$ can be naturally identified with
the hyperbolic space ${\mh}^n$ with linear metric induced by $h_{0}.$
Thus in order to compute $M(\zeta)$ we need only compute for
$R_{0}(\zeta).$ 
 It is well known, see for example Lemma 2.1 of \cite{gz3}, that
\begin{gather}
\begin{gathered}
R_0(\zeta,x,x',y,y')=
\left(\frac{1}{2} \pi^{-\frac{n}{2}}
\frac{\Gamma(\zeta)}{\Gamma\left(\zeta-\frac{n-2}{2}\right)}\right)
\frac{x^\zeta{x'}^{\zeta}}{(x^2+{x'}^2+|y-y'|_0^2)^{\zeta}} 
\left||h_0|\frac{dx}{x}\frac{dy}{x^n}\frac{dx'}{x'}\frac{dy'}{{x'}^n} \right|^{\ha}
+S_1 \\ S_1 \in {\mathcal{A}}^{\zeta+1, \zeta+1},
\end{gathered}\label{princ1}
\end{gather}
where $|y-y'|_0$ is the distance in the $h_0$ metric and 
$|h_0|$ denotes its volume element. (Here we have multiplied by the
appropriate half-density.) 
Since $\rho=x(x^2+{x'}^2+|y-y'|^2)^{-\ha}$ and 
$\rho'=x(x^2+{x'}^2+|y-y'|^2)^{-\ha}$ we deduce from \eqref{princ1} that
$R_0(\zeta)=\rho^\zeta{\rho'}^\zeta \mu$ where $\mu$ is the half density 
induced on
the front face.  By an abuse of notation we denote the restrictions of
$\rho$ and $\rho'$ to the front face also by $\rho$ and $\rho'.$
Thus $F|_{\mcf}$ is just the half-density induced  on
$\mcf.$ 
This concludes the proof of the Proposition.
\end{pf}

Now it follows from \eqref{POgb} and Proposition \ref{PF3} 
that $x^{-\zeta+\frac{n}{2}}E(\zeta),$ with $2\zeta \not\in \mz,$ and $\zeta$
not a pole of $R(\zeta),$ is a smooth section of 
$x^{-\ha}\Gamma^{\ha}(X \times \p X).$ Since
$x^{-\frac{n+1}{2}} \Gamma^{\ha}(X \times \p X)=
\Gamma_0^{\ha}(X \times \p X),$ we have that
$x^{-2\zeta}E$ is a smooth section of $\Gamma_0^{\ha}(X \times \p X).$
Therefore we have from \eqref{PF4} that
\begin{cor}\label{PO5} For $2\zeta \not\in \mz$ and $\zeta$ not a pole of 
$R(\zeta),$ the Eisenstein function, $E(\zeta),$ defined by 
\eqref{PO2}, is a smooth section of
$\Gamma_0^{\ha}(X\times \p X),$  
which is holomorphic in $\zeta.$  Moreover, for any product decomposition
$X \sim \p X \times [0,\eps),$ and for any
$f,g \in C^\infty(\p X, \Gamma^{\ha}(\p X)),$ we have that,
 as $x \downarrow 0,$
\begin{gather}
\lan E,f\otimes g\ran (x)= 
x^{-\frac{n}{2}}\left(x^{\zeta}h_{1,\zeta}(x) + 
x^{n-\zeta}h_{2,\zeta}(x)\right)
\left|\frac{dx}{x}\right|^{\ha}, \;\ x>0
\label{PO6}
\end{gather}
where $h_{i,\zeta}\in C^\infty([0,\eps)),$  $i=1,2,$ depend holomorphically in
$\zeta.$
\end{cor}

We observe that,
as an element of $C^\infty(X\times \p X, \Gamma_0^{\ha}(X\times \p X)),$
$E(\zeta)$ defines, by duality, a map  
\begin{gather*}
E(\zeta): C^\infty(\p X, \Gamma^{\ha}(\p X)) \longrightarrow
C^\infty(X, \Gamma_0^{\ha}(X)) \\
\lan E(\zeta)(f), v\ran = \lan E, f \otimes v\ran, \;\ 
f \in C^\infty(\p X, \Gamma^{\ha}(\p X)), \;\
v \in C^\infty(X, \Gamma_0^{\ha}(X)).
\end{gather*}
By definition of the resolvent, $R(\zeta),$ the kernel of
$\left(\laphd_x+V(x) +\zeta(\zeta-n)\right) R(\zeta)$ is supported on the 
diagonal in $X \times X.$ In particular we find from the definition of
$E(\zeta),$ that if 
$f \in C^\infty(\p X, \Gamma^{\ha}(\p X)),$
\begin{gather*}
\left(\Delta_x+V(x) +\zeta(\zeta-n)\right) \left(E(\zeta)f\right)=0 \text{ in }
X.
\end{gather*}

Moreover it follows from \eqref{PF5} and \eqref{PO6} that for
any  $f \in C^\infty(\p X, \Gamma^{\ha}(\p X)),$
\begin{gather*}
\left(E(\zeta)f\right)(x,\cdot)= x^{\zeta}f_{+}+ x^{n-\zeta} f_{-}, \;\
f_{\pm} \in C^{\infty}(X, \Gamma_{0}^{\ha}(X)), \;\ 
x^{\frac{n}{2}}f_{-}|_{\p X}=M(\zeta)f,
\end{gather*}
where $M(\zeta),$  is given by Proposition \ref{G01}.
This shows that $\frac{1}{M(\zeta)}E(\zeta)$  is the Schwartz kernel of the 
Poisson operator. 

For completeness, as the general result does not seem to be 
in the literature, we prove the uniqueness of the generalized eigenfunction
$E(\zeta)f$. 
The case $\Re\zeta=\frac{n}{2}$ has been proved by Borthwick in 
\cite{Bo}.  Our proof, which is based on an argument of
\cite{masy}, is not very different from his.
\begin{prop} \label{prop:eiguniqueness}
Let $\zeta \in \mc$ be such that $2\zeta\not\in \mz, $ 
$\zeta \not \in (-\infty,\frac{n}{2}],$ 
and $\zeta(\zeta-n)$ is not
 in the point spectrum of $\Delta_g.$  
 Suppose that 
$u = x^{\zeta}f+x^{n-\zeta}f'$ with $f, f'
\in C^{\infty}(X),$  satisfies
 $(\Delta + \zeta(\zeta-n))u=0.$  If  $f'|_{\p X}=0,$ then $u=0.$
\end{prop}
\begin{pf}
Substituting $u = x^{\zeta}f+x^{n-\zeta}f'$ in the equation 
$(\Delta - \zeta(n-\zeta))u=0,$  equating the powers of $x,$ and using that
$2\zeta\not\in \mz,$ we deduce that
if
$f'|_{\p X}=0$ then, in fact $f'$ vanishes to infinite order at $\p
X,$ and so can be absorbed into $f.$ 
 So we may assume that  $u=x^\zeta f.$

If  $\Re \zeta > n/2$ then $ u$ is an $L^{2}$ eigenfunction and,
by our assumption on $\zeta,$ must be zero. 

To analyze the case  $\Re \zeta \leq n/2,$ we proceed as in \cite{masy}. Let
$\phi \in C^{\infty}(\mr),$ $\phi(t)\geq 0,$ $\phi'(t) \geq 0,$ with
 $\phi(t)=0$ for  $t<1 $ and $\phi(t)=1$ for $t>2,$ and let
$(x,y)$ define a product decomposition near the boundary as in Proposition
\ref{prop:modelform}. Then $\phi(\eps^{-1}x) u \in C^\infty(X)$ vanishes near
$\p X$ and
the self-adjointness of $\Delta_g$ gives that
\begin{gather*}
\int \limits_{X} \left(\left[ \Delta_g ,\phi(\eps^{-1} x)\right] u\right)
\bar{u}\; dg= 
\int \limits_{X} \left( ( \Delta_g \phi(\eps^{-1} x)u)\bar{u} - 
(\phi(\eps^{-1} x) \Delta_g u)\bar{u} \right)\; dg
=\\ 2i\; \Im\;[\zeta(\zeta-n)]\int\limits_{X}\phi(\eps^{-1} x)|u|^2 \; dg ,
\end{gather*}
where $dg$ is the Riemannian measure induced by the density.

Since in this product decomposition
$$\Delta_{g}=-(x\p_x)^2+nx\p_x -x^2F(x,y)\p_x+ x^{2}Q\left(x,y,\p_{y}\right),
$$
with $F$ smooth, we obtain,
\begin{gather*}
\int \limits_{X} \left(\left[ \Delta_g ,\phi(\eps^{-1} x)\right] u\right)
\bar{u}\; dg= \int \limits_{ X}\left([ -(x\p_x)^2 +nx\p_x-x^2F(x,y)\p_x, \phi(x/\eps) ] u 
\right)\bar{u}\; \frac{dx}{x^{n+1}}\; dh= \\
2i\; \Im\; [\zeta(\zeta-n)]\int\limits_{X}\phi(\eps^{-1} x)|u|^2
\; dg,
\end{gather*} 
where $h$ is the natural density induced by $g$ in $y.$

Now if we have $u = x^{\zeta} f,$  then, after setting $x=\eps\tau,$ 
integrating by parts, and using that $\phi(2)=1$ and $\phi(1)=0,$  we 
obtain
\begin{gather}
\begin{gathered}
\left(-2i (2\eps)^{2\Re\zeta-n}\; \Im\zeta 
+2i \eps^{2\Re\zeta-n}\; \Im\zeta\; (2\Re\zeta-n)
 \int_{1}^{ 2} \tau^{2\Re\zeta-n-1}\phi(\tau) \; d\tau\right)
\int\limits_{\p X}|f|^{2}(0,y)\; dh +O\left(\eps^{2\Re\zeta-n+1}\right)
= \\ 2i\; \Im\; [\zeta(\zeta-n)]\int\limits_{X}\phi(\eps^{-1} x)|u|^2 \; dg.
\end{gathered}\label{uniq0}
\end{gather}

Observe that 
\begin{gather}
\begin{gathered}
\int\limits_{X}\phi(\eps^{-1} x)|u|^2 \; dg=
\int_{\eps}^{2\eps}\int_{\p X}\phi(\eps^{-1} x)x^{2\Re\zeta-n-1}|f|^2(x,y)\; 
dx dh + O(1)= \\
\eps^{2\Re\zeta-n}\int_{1}^{2} \phi(\tau) \tau^{2\Re\zeta-n-1}\; d\tau
\int_{\p X} |f|^2(0,y)dh + O\left(\eps^{2\Re\zeta-n+1}\right) + O(1).
\end{gathered}\label{uniq2}
\end{gather}

Since $\Im\zeta(2\Re\zeta-n)=\Im\;[\zeta(\zeta-n)],$ we deduce from 
\eqref{uniq0}
and \eqref{uniq2} that
\begin{gather}
-2i(2\eps)^{2\Re\zeta-n} \;\Im\zeta\;
\int \limits_{\p X}
|f|^{2}(0,y)\; dh+O\left(\eps^{2\Re\zeta-n+1}\right)=
2i\; \Im\; [\zeta(\zeta-n)] O(1). \label{uniq1}
\end{gather}
When $2\Re\zeta<n,$ since this holds as 
$\eps\rightarrow 0,$ we deduce that $f|_{\p X}=0.$
Observe that when $2\Re\zeta=n,$ $\Im\; [\zeta(\zeta-n)]=0.$ Thus the right 
hand side of \eqref{uniq1} vanishes.  Letting $\eps \rightarrow 0,$ we also
deduce that $f|_{\p X}=0.$

Once $f|_{\p X}$ vanishes it follows, using the indicial equation, and the 
fact that $2\zeta \not \in \mz,$ that $u$ must vanish to infinite
order at the boundary. Thus that $u \in L^2(X)$ and therefore $u=0.$
\end{pf}

The scattering matrix, acting on half-densities, can then be defined, for
the values of $\zeta$ as in Proposition \ref{prop:eiguniqueness}, and
such that $M(\zeta)\not=0,$ as the map
\begin{gather*}
S(\zeta): \Gamma^{\ha}(\p X) \longrightarrow \Gamma^{\ha}(\p X) \\
S(\zeta)f =\frac{1}{M(\zeta)} x^{\frac{n}{2}} {f_{+}}|_{\p X},
\end{gather*}
with $M(\zeta)$ defined as above.
Thus it follows from the first equation in \eqref{PF5} that
\begin{prop}\label{PO7} For the values of $\zeta$ as in Proposition 
\ref{prop:eiguniqueness}, and such that $M(\zeta)\not=0,$
the scattering matrix $S(\zeta)$ is a 
pseudo-differential operator in $\p X,$ acting on half-densities,
which is meromorphic in $\zeta.$
Moreover its kernel, which we also denote by $S(\zeta),$ satisfies
\begin{gather}
\beta_{\p}^{*} S(\zeta)= \frac{1}{M(\zeta)}
\beta^*\left( x^{-\zeta+\frac{n}{2}} 
{x'}^{-\zeta+\frac{n}{2}} R(\zeta)\right)|_{T \cap B},\label{PO8}
\end{gather}
where $T \cap B$ is the intersection of the top and bottom faces, and
$M(\zeta)$ is defined in Proposition \ref{G01}.
\end{prop}

Observe that  the right hand side of \eqref{PO8} gives a meromorphic
 extension of $S(\zeta)$ for values of $\zeta$ that are not poles of
$\frac{1}{M(\zeta)} R(\zeta).$

As pointed out in the introduction, this definition of the scattering matrix is
dependent on the choice of the defining function $x.$  There is a standard way
to remove, see for example \cite{gz1,per}, and view it as an operator 
\begin{gather*}
S(\zeta) : 
C^{\infty}(\p X, \Gamma^{\ha}(\p X) \otimes |N^{*}(\p X)|^{n-\zeta}) 
\longrightarrow
  C^{\infty}(\p X, \Gamma^{\ha}(\p X) \otimes |N^{*}(\p X)|^{\zeta}).
\end{gather*}

Now this scattering matrix is not quite the same as the one 
defined in the introduction, as this is the scattering matrix
associated to the operator acting on half-densities rather than on
functions.  Let $\omega_{0}$ denote the canonical density over the
boundary induced by $h.$ 
To get the appropriate Eisenstein function for functions
we take, $\omega^{-\ha}(x,y) E(\zeta) \omega_{0}^{\ha}(y').$ We thus
see that the scattering matrix on functions is obtained by trivializing
the half-density bundle over the boundary by $\omega_{0}^{\ha}.$ 
Note that conjugating the scattering matrix by the trivializing
half-density will not affect the principal symbol nor it will affect
the principal symbol of the difference of two scattering matrices associated
to differing metrics which agree at the boundary so in the next
section where we establish our inverse result it is irrelevant which
definition we use.

\section{The Principal Symbol}\label{thps}

We compute the principal symbols of $S(\zeta)$ and $S_1(\zeta)-S_2(\zeta).$
Throughout this section we 
assume that $\zeta$ is not a pole of the right hand side of \eqref{PO8}. We 
also fix a product structure in which 
\begin{gather}
\begin{gathered}
g_{j}=\frac{dx^2 +h_j(x,y,dy)}{x^2}, \;\ V_j \in C^\infty(X), \;\ j=1,2, \;\ 
h_1(x,y,dy)-h_2(x,y,dy)=x^k L(y,dy)+O(x^{k+1}),\\ V_j(0,y)=0, \;\ j=1,2, \;\;
V_1-V_2=x^kW(y)+O(x^{k+1}) .
\end{gathered}\label{hyp}
\end{gather}

First, we prove Theorem \ref{pssm0}.
\begin{pf} It follows from \eqref{PO8} and \eqref{PFA9} that the leading 
singularity of $\beta_{\p}^*S(\zeta)$ is given by 
$\frac{1}{M(\zeta)}F|_{T\cap B} R^{-2\zeta}.$ As observed in the proof of 
Proposition \ref{G01}, $F|_{T\cap B}$ is the induced half-density on
$T\cap B.$ Thus, pushing forward to $\p X \times \p X,$ gives that the leading
singularity of $S(\zeta)$ is given by
$\frac{1}{M(\zeta)}|y-y'|^{-2\zeta}$ times 
half-density given by $h.$ The density term in $M(\zeta)$ cancels with
that of $h.$ Taking the Fourier transform we find that
the principal symbol of $S(\zeta)$ is given by $C(\zeta)|\xi|^{2\zeta-n},$
where $|\xi|$ is the length of the covector $\xi$ with respect to the metric 
induced by $h.$  Note that the
principal symbol could also be computed by observing that it must
agree with that in the almost product case and that doing so gives the
explicit value of the constant - we have proceeded in the other way in
order to prepare the ground for our next result. 
\end{pf}

As a consequence of Theorem \ref{dres} and Proposition \ref{PO7} we obtain
\begin{prop}\label{PS1} Let $g_j,$ $V_j,$ $j=1,2,$ satisfy \eqref{hyp}.
Let $S_j(\zeta),$ $j=1,2$ be the scattering matrix 
corresponding to $g_j, V_j.$ Let $M(\zeta)$ be defined as above. Then
\begin{gather*}
S_2(\zeta)-S_1(\zeta)= \frac{1}{M(\zeta)}
\left(\Lambda_1(\zeta)+\Lambda_2(\zeta)\right),
\end{gather*}
where
$\Lambda_2 \in \Psi^{-\infty}(\p X, \Gamma^{\ha}(\p X)),$ and the 
Schwartz kernel of $\Lambda_1$ satisfies
\begin{gather*}
\beta_{\p}^{*} \Lambda_1(\zeta)= 
\left(R^{k-2\zeta+n}\rho^{\frac{n}{2}}
{\rho'}^{\frac{n}{2}} \alpha(\zeta)\right)|_{\rho=\rho'=0},
\end{gather*}
with $\alpha(\zeta)$ is defined by \eqref{dres1}.
\end{prop}
\begin{pf}  We will apply \eqref{dres01} and \eqref{dres1}
 to \eqref{PO8}.
Since the lift of the Schwartz kernel of $G_2$ defined in
\eqref{dres01} and \eqref{dres1}, under $\beta$ 
vanishes to infinite order at the top and bottom faces,  it does not 
contribute to the difference of the scattering matrices. Also
notice that if 
$\gamma \in C^{\infty}(X\times X, \Gamma_0^{\ha}(X\times X)),$ then
\begin{gather*}
\left(x^{-\zeta+\frac{n}{2}}{x'}^{-\zeta+\frac{n}{2}} 
x^{\zeta}{x'}^{\zeta}\gamma\right)|_{x=x'=0} \in 
C^{\infty}(\p X \times \p X, \Gamma^{\ha}(\p X \times \p X)).
\end{gather*}
So $G_3(\zeta)$ contributes to the difference of the scattering matrices with
a smoothing operator.  
Finally observe that 
\begin{gather}
\beta^*( x^{-\zeta+\frac{n}{2}}{x'}^{-\zeta+\frac{n}{2}} G_1)|_{T\cap B}=
\left(R^{k-2\zeta+n}\rho^{\frac{n}{2}}
{\rho'}^{\frac{n}{2}} \alpha(\zeta)\right)|_{\rho=\rho'=0}. \label{PS2}
\end{gather}
This concludes the proof of the proposition.
\end{pf}

Next we compute the leading singularity of $S_2(\zeta)-S_1(\zeta).$
The main part of the calculation is

\begin{lemma}\label{PS4}  Let $g_j,$  $V_j,$ $j=1,2,$ satisfy \eqref{hyp} and
let $S_j,$ $j=1,2$ be the scattering matrix corresponding to $g_j,$ $V_j.$ Let
$p \in \p X$ and assume that, after a linear
transformation, $h_0(p)=\operatorname{Id}.$ Let $S_j(\zeta),$ $j=1,2,$ be 
the scattering matrices acting on half-densities.  Then, for $M(\zeta)$ as 
above, 
\begin{gather*}
S_2(\zeta)-S_1(\zeta)= \frac{1}{M(\zeta)}\left(B_1(\zeta)+B_2(\zeta)\right),
\end{gather*}
where in local coordinates $x,y',$ valid near $p=y,$ with $w=y-y',$
$\rho=x/|w|,$ $\rho'=x'/|w|,$ $Y=w/|w|=(y-y')/|y-y'|,$
valid near $T \cap B,$ the lift of the kernels of $B_1$ and $B_2$ under 
$\beta_{\p}$
are given by
\begin{gather}
\begin{gathered}
\beta_{\p}^{*}B_1= |w|^{k-2\zeta+n}\alpha(\zeta,0,\frac{w}{|w|},y,0,0)
\left|\frac{d|w|}{|w|}\frac{dY}{|w|^n} dy'\right|^{\ha}, \\
\beta_{\p}^{*}B_2= |w|^{k-2\zeta+n+1}
\widetilde{\alpha}(\zeta,|w|,\frac{w}{|w|},y,0,0)
\left|\frac{d|w|}{|w|}\frac{dY}{|w|^n}dy'\right|^{\ha}, \;\ 
\widetilde{\alpha} \text{ smooth}.
\end{gathered}\label{TS1}
\end{gather}
Moreover, for $2\Re\zeta\geq \text{max}(k+2, n-k+1),$ 
\begin{gather}
\begin{gathered}
\alpha(\zeta,0,\frac{w}{|w|},y,0,0)= \\
C(\zeta)\left[T_1(k,\zeta)
\sum_{i,j=1}^{n}H_{ij}(y)|w|^{2\zeta-k}\p_{w_i}\p_{w_j}|w|^{k+2-2\zeta} +
T_2(k,\zeta)\left(W(y)-\oq k(n-k)T(y)\right) \right] \\
C(\zeta)=\left(\ha \pi^{-\frac{n}{2}} 
\frac{\Gamma(\zeta)}{\Gamma\left(\zeta-\frac{n-2}{2}\right)}\right)^2, \\
T_l(k,\zeta)= \int_{0}^{\infty}\int_{\mrn}
\frac{u^{2\zeta+k+3-2l-n}}{(u^2+|V|^2)^{\zeta}
(u^2+|e_1-V|^2)^{\zeta}} dV d u, \;\ e_1=(1,0,...,0), \;\ l=1,2.
\end{gathered} \label{sm1}
\end{gather}
\end{lemma}
\begin{pf} In these coordinates, \eqref{PS2} is given by
\begin{gather}
\beta^*( x^{\zeta+\frac{n}{2}}{x'}^{\zeta+\frac{n}{2}} G_1)|_{T\cap B}=
|w|^{k-2\zeta+n}
\alpha(\zeta,|w|,Y,y,0,0)\left|\frac{d|w|}{|w|}
 \frac{dY}{|w|^n} dy' \right|^{\ha}.\label{PS3}
\end{gather}
 Now we use Proposition \ref{PS1} and observe that
$\beta_{\p}= \beta|_{T \cap B}.$  Equation \eqref{TS1} is just
the first order Taylor's expansion in $|w|$ of the function 
$\alpha(\zeta,|w|,Y,0,0).$

We observe that
$\alpha(\zeta,0,Y,y,0,0)|dY dy'|^{\ha}$ 
is the restriction of $R^{\frac{n}{2}}\rho^{\frac{n}{2}}{\rho'}^{\frac{n}{2}}
\alpha(\zeta)$ to the the intersection of the top, bottom and front faces, 
$T\cap B \cap F=\{R=\rho=\rho'=0\}.$   
We know from Theorem \ref{dres}  
that the half-density $R^{\frac{n}{2}}\alpha(\zeta),$ restricted to the front
face, satisfies to \eqref{dres02}. By Remark \ref{remark1} this equation has 
a unique solution, and it can be solved directly.  Then we find the
value of a solution to \eqref{dres02} at $T\cap B \cap F.$
Instead of coordinates $R,$ $Y$ and $\rho,$ it is convenient to use
$s=\frac{x}{x'}, $ $z=\frac{y-y'}{x'}.$  Then
the front face is given by $x'=0$ and we have
\begin{gather}
\begin{gathered}
\rho=\frac{x}{(x^2+{x'}^2+|y-y'|^2)^{\ha}}=\frac{s}{(1+s^2+|z|^2)^{\ha}}, \\
\rho'=\frac{x'}{(x^2+{x'}^2+|y-y'|^2)^{\ha}}=\frac{1}{(1+s^2+|z|^2)^{\ha}}.
\end{gathered}\label{PS51}
\end{gather}
The intersection of the top, bottom and front faces, $T\cap B \cap F,$ is then
 given by $\{x'=s=0, \; |z|=\infty \},$ and since $h_0=\operatorname{Id},$
 equation \eqref{dres02} is reduced to
\begin{gather}
(\Delta + \zeta(\zeta-n))
\left(s^{\zeta}(1+s^2+|z|^2)^{\frac{k-2\zeta}{2}} \alpha(s,z)\right)= 
N_p(s^k E)G, \label{PS6}
\end{gather}
where $\Delta$ is the Laplacian in the hyperbolic space. Hence we have
\begin{gather} 
s^{\zeta}(1+s^2+|z|^2)^{\frac{k-2\zeta}{2} }\alpha(s,z)=
G(N_p(s^k E)G)(s,z). \label{PS.new16}
\end{gather}
We recall that the uniqueness of the solution to \eqref{PS.new16} is 
established in  Remark \ref{remark1}.
It follows from \eqref{AR.2} that
\begin{gather}
N_p(E)=\sum_{i,j=1}^{n} H_{ij}(y) s\p_{z_i} s\p_{z_j} + 
\left(W(y)-\frac{1}{4}k(n-k)T(y)\right). \label{PS8}
\end{gather}
We recall from Lemma 2.1 of \cite{gz3} that 
\begin{gather}
G(s,z)= \left(\frac{1}{2} \pi^{-\frac{n}{2}}
\frac{\Gamma(\zeta)}{\Gamma\left(\zeta-\frac{n-2}{2}\right)} 
\frac{s^\zeta}{(1+s^2+|z|^2)^{\zeta}} \right)\left|\frac{ds}{s} 
\frac{dz}{s^n} dy'\right|^{\ha} +G_1 \label{PS7} 
\end{gather}
where $G_1$ has a conormal singularity at $\{s=1, z=0\}$ and,
 near the boundary, $G_1 \in {\mathcal{A}}^{\zeta+1, \zeta+1},$
where ${\mathcal{A}}^{a,b}$ denotes the space of
half-densities of the form $s^a |z|^{-b}|\frac{ds}{s} 
\frac{dz}{s^n} dy'|^{\ha}.$
It follows from Proposition 6.19 of \cite{mm} that
$G\left({\mathcal{A}}^{\zeta+k+1, \zeta+1}\right) \subset 
{\mathcal{A}}^{\zeta, \zeta-k+1}.$
Since, as in equation (4.12) of \cite{mm}, $G$ acts as a convolution operator 
with respect to the group action defined in section 3 of that paper, we find 
that, 
\begin{gather}
\begin{gathered} 
s^{\zeta}(1+s^2+|z|^2)^{\frac{k-2\zeta}{2}} \alpha(s,z)\left|\frac{ds}{s} 
\frac{dz}{s^n}dy'\right|^{\ha}= \\
C(\zeta)\left( \sum_{i,j=1}^{n} H_{ij}(y) \p_{z_i} \p_{z_j} I_1(k,\zeta,s,z) +
\left( W(y)-\oq k(n-k)T(y)\right) I_2(k,\zeta,s,z) 
 \right) \left|\frac{ds}{s} \frac{dz}{s^n} dy' \right|^{\ha} +\beta, \\
\text{ where } C(\zeta)=\left(\frac{1}{2} \pi^{-\frac{n}{2}}
\frac{\Gamma(\zeta)}{\Gamma\left(\zeta-\frac{n-2}{2}\right)}\right)^2, \;\
 \beta \in {\mathcal{A}}^{\zeta, \zeta-k+1}, \text{ and } \\
I_{l}(k,\zeta,s,z)= \int_0^{\infty}\int_{\mrn} 
\frac{t^{\zeta}}{(1+t^2+|U|^2)^{\zeta}(1+\frac{s^2}{t^2} +
|z-\frac{s}{t}U|^2)^{\zeta}}\left(\frac{s}{t}\right)^{\zeta+k+4-2l}
 \frac{dt}{t} dU.
\end{gathered}\label{PS9}
\end{gather}
Recall that our goal is to compute the restriction of 
$R^{\frac{n}{2}}\rho^{\frac{n}{2}}{\rho'}^{\frac{n}{2}} \alpha(\zeta)$ to
$T\cap B \cap F.$ In these coordinates
\begin{gather*}
R^{\frac{n}{2}}\rho^{\frac{n}{2}}{\rho'}^{\frac{n}{2}}= 
{x'}^{\frac{n}{2}}\frac{s^{\frac{n}{2}}}{(1+s^2+|z|^2)^{\frac{n}{4}}}.
\end{gather*}
So, after restricting to the front face, which is given by $\{x'=0\},$
we have to restrict 
$$\frac{s^{\frac{n}{2}}}{(1+s^2+|z|^2)^{\frac{n}{4}}}\alpha(s,z)\left|\frac{ds}{s}
\frac{dz}{s^n} dy'\right|^{\ha}$$ to the corner $T\cap B \cap F=\{s=0, |z|=\infty\}.$
This is the same as the restriction of
$$\frac{s^{\frac{n}{2}}}{|z|^{\frac{n}{2}}}\alpha(s,z)\left|\frac{ds}{s}
\frac{dz}{s^n} dy'\right|^{\ha}.$$ 
Notice that for $Y=z/|z|,$ we have $dY=dz/|z|^n.$
Thus the
value of $s^{\frac{n}{2}}|z|^{-\frac{n}{2}}\alpha$ at 
$T\cap B \cap F$ is then given by $A(Y)|dY dy'|^{\ha},$ $Y=z/|z|,$ where
\begin{gather}
\begin{gathered}
A(Y)= \\ \lim \limits_{s\rightarrow 0, |z| \rightarrow \infty}
\frac{1}{s^{\zeta}(1+s^2+|z|^2)^{\frac{k-2\zeta}{2}}}
\left[\sum_{i,j=1}^{n} H_{ij}(y) \p_{z_i} \p_{z_j} I_1 + 
\left(W(y) -\frac{1}{4}k(n-k)T(y)\right) I_2\right] .
\end{gathered}\label{PS10}
\end{gather}

Set $|z|u=\frac{s}{t}$ and $U=\frac{t}{s}|z|V$ and observe that 
$I_l(k,\zeta,s,z)=I_l(k,\zeta,s,|z|),$ so we can also 
set $z=|z|e_1,$ $e_1=(1,0,...,0).$  Then
\begin{gather*}
I_l(k,\zeta,s,z)=s^{\zeta}|z|^{-2\zeta+k+4-2l}\int_{0}^{\infty}\int_{\mrn}
\frac{u^{2\zeta+k+3-2l-n}}{(u^2+\frac{s^2}{|z|^2}+|V|^2)^{\zeta}
(\frac{1}{|z|^2}+u^2+(e_1-V)^2)^{\zeta}} dV du.
\end{gather*}
To analyze the limit of $I_l(k,\zeta,s,z)$ as $s\rightarrow 0$ and
$|z|\rightarrow \infty,$ we begin by proving
\begin{lemma}\label{limit}
For $k\geq 1,$ and for 
$2\Re\zeta\geq \text{max}\left(n-k+1,k+2\right),$ we have
\begin{gather*}
J(l,k,\zeta)=\int_{0}^{\infty}\int_{\mrn}
\frac{u^{2\Re\zeta+k+3-2l-n}}{(u^2+|V|^2)^{\Re \zeta}
(u^2+(e_1-V)^2)^{\Re\zeta}} dV d u <\infty.
\end{gather*}
\end{lemma}
\begin{pf} Observe that for $V=(v, V'),$ $V'\in \mr^{n-1}$
and $|V'|=\rho,$
\begin{gather*}
J(l,k,\zeta))=
|\ms^{n-2}| \int_{0}^{\infty} \int_{0}^{\infty} \int_{\mr}
\frac{u^{2\Re\zeta+k+3-2l-n}\rho^{n-2}}{(u^2+\rho^2+v^2)^{\Re\zeta}
(u^2+\rho^2+(v-1)^2)^{\Re\zeta}} dv du d\rho.
\end{gather*}
Setting $v=R\cos \phi,$ $u=R\sin\phi\cos\theta,$
$\rho=R\sin\phi\sin\theta,$ $0<\phi<\pi,$ $0<\theta<\frac{\pi}{2},$ we obtain
\begin{gather*}
J(l,k,\zeta)= K(\zeta)
\int_{0}^{\infty}\int_{0}^{\pi}
\frac{R^{k+3-2l}(\sin\phi)^{2\Re\zeta+k+2-2l}}{\left[(R-\cos\phi)^2+
(\sin\phi)^2\right]^{\Re\zeta}} d\phi dR, \\
K(\zeta)=|\ms^{n-2}|\int_{0}^{\frac{\pi}{2}}(\cos\theta)^{2\Re\zeta+k+3-2l}
(\sin\theta)^{n-2}d\theta .
\end{gather*}
Thus, for $k\geq 2$ and  $\zeta$ as above,  we have that
\begin{gather*}
J(l,k,\zeta) \leq K_1(\zeta)( \int_{0}^{4} \int_{0}^{\pi} R^{k+3-2l} 
(\sin\phi)^{k+2-2l} d\phi dR + \\
K_2(\zeta) \int_{4}^{\infty}\int_{0}^{\pi} R^{-2\Re\zeta+k+3-2l}
(\sin \phi)^{2\Re\zeta+k+2-2l} d\phi dR )<\infty
\end{gather*}
The same argument can be used to show that $J(1,1,\zeta)<\infty.$
When $k=1$ and $l=2,$ another argument has to be used. 
Setting
$R=\cos\phi +t \sin\phi$ we find that
\begin{gather*}
J(2,1,\zeta) \leq K(\zeta) \int_{-\infty}^{\infty}\int_{0}^{\pi}
(1+t^2)^{-\Re\zeta} d\phi dt <\infty.
\end{gather*}
This concludes the proof of the Lemma.
\end{pf}

Thus the dominated convergence theorem gives that for
$T_l(k,\zeta,s,z)=s^{-\zeta}|z|^{-2\zeta+k-4+2l}I_l(k,\zeta,s,z),$
\begin{gather}
\lim_{s\rightarrow 0, |z|\rightarrow \infty} T_{l}(s,z)= T_l(k,\zeta)=
\int_{0}^{\infty}\int_{\mrn}
\frac{u^{2\zeta+k+3-2l-n}}{(u^2+|V|^2)^{\zeta}
(u^2+(e_1-V)^2)^{\zeta}} dV d u. 
\label{PS101}
\end{gather}

By identical considerations we deduce that 
\begin{gather*}
\p_{z_j}T_l(s,z)=O\left(\frac{s^2}{|z|^3}\right), \;\
\p_{z_m}\p_{z_j}T_l(s,z)=O\left(\frac{s^2}{|z|^4}\right).
\end{gather*}
Hence
\begin{gather}
\begin{gathered}
\p_{z_i}\p_{z_j}I_1(s,z)=
C_1(\zeta)s^{\zeta}(\p_{z_i}\p_{z_j}|z|^{k+2-2\zeta})T_1(s,z)+
O\left(s^{\Re\zeta}|z|^{k-2-2\Re\zeta}\right).
\end{gathered}\label{PS11}
\end{gather}

Using that $z=\frac{y-y'}{x'}$ we find that 
$\frac{z}{|z|}=\frac{w}{|w|}=\frac{y-y'}{|y-y'|}.$ Therefore
 \eqref{sm1} follows directly form \eqref{PS9},
\eqref{PS10}, \eqref{PS101} and 
\eqref{PS11}. This concludes the proof of the Proposition.
\end{pf}

Now we can prove Theorem \ref{pssm}. 
\begin{pf} It follows from \eqref{TS1} and \eqref{sm1}, 
the leading singularity of the kernel of the difference of the
scattering matrices is given by
\begin{gather}
\frac{C(\zeta)}{M(\zeta)}\left(T_1(k,\zeta)\sum_{i,j=1}^{n}H_{ij}(y)
\p_{w_i}\p_{w_j}|w|^{k+2-2\zeta} +
T_2(k,\zeta)\left( W(y)-\oq k(n-k)T(y)\right)|w|^{k-2\zeta}\right)\label{last}
\end{gather}
times a non-vanishing 
smooth half-density, where $C(\zeta)$ is given by \eqref{sm1} and
$M(\zeta)$ by Proposition \ref{G01}. 
We obtain \eqref{psf} by taking Fourier transform in $w$ of \eqref{last},
and observing that \eqref{sm1} was
obtained under the assumption that $h_0=\operatorname{Id},$
and using the fact that $h_0$ is symmetric. The coefficients of $T_j(k,\zeta),$
$j=1,2$ in \eqref{psfnew} arise when we take the Fourier transform of the
corresponding power of $|w|.$ See for example page 363 of \cite{gelfand}.
 This ends the proof of the theorem.
\end{pf}
We now prove Corollaries \ref{cor1} and \ref{cor2}. The proof of Corollary
\ref{cor2} is a direct consequence of the fact that, for every $k,$
 $A_2(k,\zeta)\not=0$ for at least one value of $\zeta.$ The proof of
Corollary \ref{cor1} requires a more delicate analysis due to the presence of
the term involving $T(y).$
\begin{pf} As we are working modulo diffeomorphism invariance we can
take a product decomposition such that each $g_j$ is of the form
\eqref{pdec}. 
 Suppose $g_1$ equals $g_2$ to order $k$
near $p$  and 
suppose that the principal 
symbol of $S_1(\zeta)-S_2(\zeta)$ of order
$2\Re\zeta-n-k$ is equal to zero at $p.$ Since $V_1=V_2$ near $p,$ we find that
$W=0.$ By a linear change of variables on the tangent space to $\p X$ at $p$
we may assume that $h_0=Id.$ 
It is clear from \eqref{psf} that if the trace is zero and
$A_{k,\zeta}$ is non-zero then
$L_{ij}(p)=0$ is zero so we need only show that off a discrete set
these hold.  
By taking $\xi=e_j=(0,...,0,1,0,...,0),$ $1$ in 
the $j$-th entry, we deduce from \eqref{psf} that
\begin{gather*}
A_1(k,\zeta)L_{ij}(p)-\frac{1}{4}k(n-k)A_2(k,\zeta)T(p)=0, \;\ 
1\leq i,j \leq n.
\end{gather*}
By taking $i=j$ and adding in $j$ we obtain, for all $\zeta$ which
is not a pole of $A_j(k,\zeta),$ $j=1,2,$
\begin{gather*}
\left(A_1(k,\zeta)-\frac{1}{4}nk(n-k)A_2(k,\zeta)\right)T(p)=0, 
\end{gather*}
Using the formulas for $A_1$ and $A_2$ given by \eqref{psfnew} and the
fact that $\Gamma(z+1)=z\Gamma(z)$ we have, again
for all $\zeta$ which
is not a pole of $A_j(k,\zeta),$ $j=1,2,$
\begin{gather*}
\left(T_1(k,\zeta)(k+2-2\zeta)(k-2\zeta+n)-
\frac{1}{4}nk(n-k)T_2(k,\zeta)\right)T(p)=0.
\end{gather*}
We know from Lemma \ref{limit} that for
$k\geq 1,$ and for $2\Re\zeta \geq \text{max}\left(n-k+1,k+2\right),$
$T_1(k,\zeta)$ and $T_2(k,\zeta)$ are finite.  In particular they are finite
for $2\zeta=k+n,$ as long as $n\geq 2.$ It is clear from the definition of
$T_j$ that
for $2\zeta=k+n,$ $T_j(k,\zeta)>0,$ $j=1,2.$ Hence $T(p)=0$ and $H_{ij}=0.$

For $n=1$ we have that,
since $k\geq 1,$  $-nk(n-k)\geq 0.$ On the other hand, for
$\zeta$ large and real $(k+2-2\zeta)(k-2\zeta+n)>0.$  Thus we also have 
$T(p)=0.$ This ends the proof of the Corollary.
\end{pf}

\section{Almost Product Type Metrics} \label{sec:prod}

In this section, we examine the scattering matrix for metrics which 
take the form,
\begin{equation}
g = \frac{dx^{2} + h(y,dy)}{x^2} + O(x^{\infty}),
\end{equation}
for some product decomposition. Our approach is analogous to that of
Christiansen, \cite{chris}, and Parnovksi, \cite{parn},  in the asymptotically Euclidean setting. 
The computation is also closely related to that of Hislop, 
\cite{His} section 2.3, for ${\Bbb H}^{n}.$ 

As we have shown in previous sections that if two metrics agree to
infinite order then the associated scattering matrices differ by a
smoothing operator, it is sufficient to compute for the manifold, 
$ \mr_{+} \times \p X,$ with metric $\frac{dx^{2} + h(y,dy)}{x^2}.$
The Laplacian is then,
$$ - \left( x \frac{\p}{\p x} \right)^{2} + n x \frac{\p}{\p x} + x^2
\Delta_{\p X},$$
where $\Delta_{\p X}$ is the Laplacian associated to $h$ on $\p X.$ 
Let $\psi_j$ be a complete orthonormal basis of eigenfunctions 
for $\Delta_{\p X}$ with
$\psi_j$ of eigenvalue of $\lambda^{2}_{j}.$ 

We then look for solutions of $(\Delta + \zeta(\zeta-n))u=0$ of the
form $x^{n/2} a(x) \psi_j (y).$ Computing as in \cite{His} we deduce
that $a$ satisfies,
$$ \left(  x^2 \frac{\p^2}{\p x^2} + x \frac{\p}{\p x} - [ x^2
\lambda^{2}_{j} + (\zeta-n/2)^{2}] \right) a(x) =0.$$ This is a modified Bessel
equation and taking the solutions which are regular at infinity, we
 see that $a$ has an asymptotic expansion  
as $x \to 0,$ and its lead term is
of the form,
$$ \frac{1}{\Gamma(1 - (\zeta-n/2))} \fracwithdelims(){ \lambda_j x }{2}^{n/2-\zeta} -
\frac{1}{\Gamma(1 + (\zeta-n/2))} \fracwithdelims(){ \lambda_j x}{2}^{\zeta-n/2}.$$ 
It now follows immediately that $S(\zeta)$ applied to $\psi_j$
multiplies it by the ratio of these coefficients:
$$- \frac{ \frac{1}{\Gamma(1 + (\zeta-n/2))} (\ha \lambda_j
)^{\zeta-n/2}}{\frac{1}{\Gamma(1 - (\zeta-n/2))} (\ha \lambda_j x
)^{n/2-\zeta}} =  2^{n-2\zeta} \frac{\Gamma(n/2 -
\zeta)}{\Gamma(\zeta-n/2)} \lambda_{j}^{2\zeta-n}.$$
As the functions $\psi_j$ form an orthonormal basis, we have now
proven the second part of 
Theorem \ref{product}.

\section{ Inverse Scattering For Black Holes}\label{schw}

We consider two models for the exterior of a static black hole,
the Schwarzschild, and the De Sitter Schwarzschild models. These are given
by
\begin{gather*}
(Y, g), \;\ Y= \mr_t \times X, \text{ where } \;\
 g=\alpha^2 dt^2 -
\alpha^{-2} dr^2 -r^2 |d\omega|^2,
\end{gather*}
$|d\omega|^2$ is the standard metric on $\ms^2.$
In the Schwarzschild model 
\begin{gather}
X=(r_+, \infty)_r \times \ms^{2}_{\omega}, \;\ \text{ and } 
\alpha=\left(1-\frac{2m}{r}\right)^{\ha}, \;\ r_+=2m<r,\label{scha}
\end{gather}
and in the De Sitter-Schwarzschild model, 
\begin{gather}
 X=(r_+, r_{++})_r \times \ms^{2}_{\omega}, \;\ \text{ and }
\alpha=\left(1-\frac{2m}{r}-\frac{1}{3}\Lambda r^2\right)^{\ha},
\;\ r_+<r<r_{++}.\label{dscha}
\end{gather}
The parameter $ m >0 $ denotes the mass of the black hole. In \eqref{dscha},
$ \Lambda,$ with $0<9m^2\Lambda<1,$ is the cosmological constant,
and $r_+,$ $r_{++}$ are the two solutions to
$\alpha=0.$

These are semi-Riemannian metrics
on the manifold with boundary $Y,$ so their Laplacians are in fact 
hyperbolic operators, we denote them $\square_g.$ We have
\begin{gather}
\square_g= \alpha^{-2}\left( D_{t}^{2}-\alpha^2r^{-2}D_r(r^2\alpha^2)D_r-
\alpha^2r^{-2}\Delta_{\omega}\right), \label{boxg}
\end{gather} 
where $D_{\bullet}=\frac{1}{i}\p_{\bullet}$ and $\Delta_{\omega}$ is the
positive Laplacian on $\ms^2.$

Therefore stationary scattering phenomena are governed by the 
operator
\begin{gather}
P=\alpha^2r^{-2}D_r(r^2\alpha^2)D_r-\alpha^2r^{-2}\Delta_{\omega}.\label{eqp}
\end{gather}

Scattering theory for the operator $P$ has been extensively studied see for 
example \cite{Ch,ChD,Ba,on,sazwsch} and the references cited there.
 It was observed in \cite{sazwsch} that, after a change of $C^\infty$ structure
on $X,$ the De Sitter-Schwarzschild model the operator $P$ can be 
viewed  as $0$-differential operator which is elliptic, and whose normal
operator is, after a linear change of variables, 
a multiple of the Laplacian on the hyperbolic space. This change in
$C^{\infty}$ structure is simply the addition of the square root of the
boundary defining function and therefore only affects smoothness up to
the boundary and not smoothness in the interior. 
Thus
the methods of \cite{mm} directly apply and it was shown in \cite{sazwsch} 
that $R(\la)=\left(P-\la^2-\frac{n^2}{4}\right)^{-1}$ has a meromorphic 
continuation to 
$\mc.$   It also follows from
the discussion in \cite{sazwsch}, and the methods of section \ref{scatmat},
that the scattering matrix can be defined in this situation.

  The case of the Schwarzschild model is more complicated. At one end, 
$\alpha=0,$ which is the black hole, the operator $P$ behaves as in
the De Sitter-Schwarzschild model, i.e, after a change in the $C^\infty$ 
structure of $X,$ it is an elliptic
$0$-differential operator and its normal operator is essentially
the hyperbolic Laplacian.
On the other end, as $r\rightarrow \infty,$ $\alpha \rightarrow 1,$ and 
the metric $g$ tends to the Lorentz metric, thus
the operator $P$ tends to the Euclidean Laplacian.  This is the case of an
asymptotically Euclidean metric. To study the scattering matrix at this
end one proceeds as in \cite{smagai}. Since the construction of the symbol of
the scattering matrix at each end only depends on the metric in a 
neighbourhood of each boundary, see \cite{smagai} and section \ref{scatmat},
it follows that modulo smoothing operators, the scattering matrices 
at each boundary are independent.

It was shown in \cite{Ba} that the resolvent $R(\la),$ for the 
Schwarzschild model,
 as an operator from ${\cal C}^\infty_0 (\stackrel{o}{X}) $ to 
${\cal C}^\infty (\stackrel{o}{X} ),$
has a meromorphic continuation from $\Im \la >0$ to ${\Bbb C} \setminus 
i \overline {\Bbb R}_- .$  It is not known whether its poles might accumulate 
at the origin.

In this section we will prove that the Taylor series of certain perturbations 
of the both models, at $\alpha=0,$ are determined from the scattering matrix 
at a fixed energy. The analogous result at $x=0$ also holds for the 
Schwarschild model, however, since its proof relies on the methods of
\cite{metric},  we will not carry it out  here.

\begin{thm}\label{pssmsch}
Let $(X, \p X)$ be a smooth manifold with boundary with dimension $n+1$,
and let $p \in \p X.$ Suppose that $g$ induces an asymptotically hyperbolic 
structures on $X$ 
and that $g = \frac{dx^2 + h(x,y,dy) }{x^2},$ with respect 
to some product decomposition near $\p X$.  
Suppose that $P$ is a smooth elliptic $0$-differential operator
of second order that its normal operator satisfies
\begin{gather}
\begin{gathered}
N_q\left(P\right)=K N_q\left(\Delta_{g}\right), \;\ \;\ 
\forall \; q \in \p X,  \\
\end{gathered} \label{sch1}
\end{gather}
where $K>0$ is a constant on each component of $\p X.$  
Then for each $\la \in \mr\setminus Q,$ $Q$ a discrete subset,
 and $f \in C^{\infty}(\p X)$ there 
exists a unique $u$ satisfying $(P-\la^2-\frac{n}{4})u=0$ of the form
\begin{gather*}
u=x^{i\la+\frac{n}{2}}f_+ + x^{-i\la+\frac{n}{2}}f_-, \;\
f_{+}|_{\p X}=f.
\end{gather*}
Moreover the scattering matrix, given by,
$$S(\la) f= f_{-}|_{\p X},$$
is a pseudo-differential operator of order $2i\lambda$.

Furthermore if  $P_2$ is another  smooth elliptic $0$-differential operator
of second order that satisfies \eqref{sch1} and is such that
\begin{gather}
P-P_2 = x^{k}\left(\sum_{i,j=1}^{n} H_{ij}(x)x\p_{y_i}x\p_{y_j} +W(x)\right)
+ O(x^{k+1}), \label{schw1n}
\end{gather}
where $H=\left(H_{ij}\right)$ is a smooth symmetric matrix.
Then
\begin{gather*}
S(\la) - S_{2}(\la) \in \Psi DO^{2i\lambda -k},
\end{gather*}
and the principal symbol of $S(\la) - S_{2}(\la)$ equals
\begin{equation}
A_1(k,\la)\sum \limits_{i,j} H_{ij}
\xi_i \xi_j |\xi|^{2i\lambda-k-2} +A_{2}(k,\la) W 
|\xi|^{2i\lambda-k}, \label{psfsch}
\end{equation}
where $h_0 = h|_{x=0},$
$|\xi|$ is the length of the co-vector $\xi$ induced by $h_0,$
and $A_1, A_2$ are  functions of $\la$ which  are not identically zero.
\end{thm}
\begin{pf} 
A line by line inspection of the proof of Theorem \ref{pssm} with
$\zeta=\frac{n}{2}+i\la$ gives the result.
\end{pf}

As an application of Theorem \ref{pssmsch} we will prove
\begin{thm}\label{invbh1} Let $X$ and $\alpha$ be given by either
\eqref{scha} or \eqref{dscha}.
Let $a_{ij}(r,\omega) \in C^{\infty}(\overline{X}),$ $0 \leq i,j \leq 2$  and
let
\begin{gather}
g=\alpha^{2} dt^{2} -\alpha^{-2}(1+\alpha a_{00}(r,\omega))dr^2 - 
\sum_{j=1}^{2} a_{0j}dr d\omega_{j} -
r^2 \sum_{i,j=1}^{2}(\delta_{ij}+\alpha a_{ij}) 
d\omega_i d\omega_j, \label{persch} 
\end{gather}
be a perturbation of the models above.
Let $X_{\ha}$ be the manifold $X$ with the new $C^{\infty}$ structure
in which $\alpha \in C^{\infty}(X_{\ha})$ is the new boundary defining 
function.  Then the operator
$P_a$ operator defined by 
\begin{gather}
\square_{g_{a}}=\alpha^{-2}\left( D_{t}^{2}- P_{a}\right)\label{eqpa}
\end{gather}
satisfies the hypotheses of Theorem \ref{pssmsch} at the boundary, 
$\{\alpha=0\},$ and there exists a product decomposition 
$(\widetilde{\alpha},\widetilde{\omega}),$ with $\widetilde{\omega}=\omega$
at $\alpha=0,$  near $\ X$ such that
for  $\la \in \mr\setminus Q,$ $Q$ a countable subset,
its scattering matrix at energy $\la$  determines the Taylor
series of $a_{ij}$ in coordinates $(\widetilde{\alpha},\widetilde{\omega})$
 at $\{\widetilde{\alpha}=0\}.$
\end{thm}
Note as before we can recover to finite order off a discrete subset
but to infinite order off a countable subset.
\begin{pf} We will only carry out the proof for the Schwarzschild model,
the other case is very similar, although the computations are more tedious,
but are essentially done in \cite{sazwsch}.

First we check the statement about the normal operator of $P_a.$ Since
$\alpha^2=1-\frac{2m}{r}$ we find that $dr=\alpha \frac{r^2}{m}d\alpha.$
Hence $g$ is given by
\begin{gather}
g=\alpha^{2} dt^2 -\frac{r^4}{m^2}(1+\alpha a_{00}(r,\omega)) d\alpha^2
- \alpha \frac{r^2}{m^2}a_{0j}d\alpha d\omega_j -
r^2 \sum_{i,j=1}^{2}(\delta_{ij}+\alpha a_{ij}) 
d\omega_i d\omega_j. \label{alpha}
\end{gather}

Let $A_0=\left(a_{ij}^{0}\right),$ where $a_{00}^{0}=\alpha^{-2},$ 
$a_{22}^{0}=a_{33}^{0}=r^2$ and $a_{ij}^{0}=0,$ $i\not = j.$
Let $A_1=\left(a_{ij}^{1}\right),$ where $a_{00}^{1}=\alpha^{-2}a_{00},$
$a_{0j}^{1}=a_{j0}^{1}=a_{j0},$ $a_{ij}^{1}=a_{ij},$
$1\leq i,j \leq 2.$   Let $A=A_0+\alpha A_1.$ Then we have
$A=A_0\left(\Id + \alpha A_{0}^{-1} A_1\right)$ and hence
\begin{gather}
\begin{gathered}
\det(A)= \det(A_0)\det\left(I+\alpha A_{0}^{-1} A_1\right)= 
\det(A_0)\left(1+\alpha T+ O(\alpha^2)\right), \;\ 
T=a_{00}+a_{11}+a_{22}, \\
A^{-1}=A_{0}^{-1}+ \alpha A_{0}^{-1} A_1 A_{0}^{-1}.
\end{gathered}\label{detf}
\end{gather}
Using \eqref{detf} and the definition of $P_a$ we find 
 that the normal operator of $P_a$ at a point $p$ at the boundary 
$\alpha=0$ is
\begin{gather*}
N_{p}\left(P_a\right)=\frac{1}{16m^2}
\left( 4(\alpha\p_{\alpha})^2 + \alpha^2\Delta_{p}\right),
\end{gather*}
where $\Delta_p$ is the Laplacian at the tangent plane to $\ms^2$ at $p.$ Thus
$N_p(P_a)$ satisfies \eqref{sch1}.

Next we consider two perturbations of the Schwarzschild metric $F$ and $H$
satisfying 
\begin{gather*}
F_{00}=\frac{r^4}{m^2}\left(1+\alpha f_{00}\right), \;\ 
F_{1j}=F_{j1}=\frac{r^{2}}{m^2}f_{1j}, \;\ 
F_{ij}=r^2\left(\delta_{ij}+\alpha r^{-2}f_{ij}\right), \\
H_{00}=\frac{r^4}{m^2}\left(1+\alpha h_{00}\right), \;\
H_{1j}=H_{j1}= \frac{r^{2}}{m^2} h_{1j}, \;\  
H_{ij}=r^{2}\left(\delta_{ij}+\alpha r^{-2}h_{ij}\right). 
\end{gather*}
Let $g_F$ and $g_H$ be defined by \eqref{alpha}, where
$f_{ij}$ and $h_{ij}$ play the r\^ole of $a_{ij}.$
Let $S_F$ and $S_H$ be the scattering matrices corresponding to $P_F$ and 
$P_H.$ 
It follows from the computation of the determinant above that, for $\alpha$ 
small,  and $f_{ij},$ $h_{ij},$ smooth,
\begin{gather*}
G_F=\frac{r^4}{m^2}(1+\alpha f_{00}(r,\omega)) d\alpha^2
+\alpha \frac{r^2}{m^2}\sum_{j=1,2} f_{0j}d\alpha d\omega_j +
r^2 \sum_{i,j=1}^{2}(\delta_{ij}+\alpha f_{ij}) 
d\omega_i d\omega_j \\
G_H=\frac{r^4}{m^2}(1+\alpha h_{00}(r,\omega)) d\alpha^2
+\alpha \frac{r^2}{m^2}\sum_{j=1,2} h_{0j}d\alpha d\omega_j +
r^2 \sum_{i,j=1}^{2}(\delta_{ij}+\alpha h_{ij}) 
d\omega_i d\omega_j 
\end{gather*}
are Riemannian metrics near $\p X.$ 

Let $(\widetilde{\alpha},\widetilde{\omega})$ be a product decomposition
of $X$ near $\p X$ in which
\begin{gather*}
G_F= d\alpha^{2} +
\widetilde{f_{ij}} d\omega_i d\omega_j, \;\
G_H= d\alpha^{2} +
\widetilde{h_{ij}} d\omega_i d\omega_j
\end{gather*}

Suppose that, in these coordinates, 
$\widetilde{f_{ij}}-\widetilde{h_{ij}}=
{\widetilde{\alpha}}^k \widetilde{u_{ij}}.$ Therefore
\begin{gather*}
P_{F}-P_{H}={\widetilde{\alpha}}^k
\left(\widetilde{u_{ij}} \p_{\omega_i} \p_{\omega_j}\right) +
O\left({\widetilde{\alpha}}^{k+1}\right).
\end{gather*}
So it follows from Theorem \ref{pssmsch} that the $k$-th order symbol
of $S_F(\la)-S_H(\la)$ determines and is determined by
$\widetilde{u_{ij}}.$

This ends the proof of the Theorem.
\end{pf}

\end{document}